\documentclass[graybox]{svmult}

\usepackage{type1cm}        % activate if the above 3 fonts are
\usepackage{makeidx}         % allows index generation
\usepackage{graphicx}        % standard LaTeX graphics tool
\graphicspath{{images/}}
\usepackage{multicol}        % used for the two-column index
\usepackage[bottom]{footmisc}% places footnotes at page bottom

\usepackage{subcaption}
\usepackage{newtxtext}       % 
\usepackage{newtxmath}       % selects Times Roman as basic font

\makeindex             % used for the subject index
\usepackage{amssymb}
\usepackage{mathtools}

\usepackage{algorithm}
\usepackage{algpseudocode}

\usepackage{tikz}
\usetikzlibrary{math}
\usetikzlibrary{spy}
\usetikzlibrary{calc}
\usetikzlibrary{patterns}
\usetikzlibrary{external}
\usepackage{hyperref}
\hypersetup{
    colorlinks=true,
    linkcolor=gray,
    filecolor=gray,      
    urlcolor=gray,
    citecolor=gray,
    pdfpagemode=FullScreen,
}
\usepackage{subfiles}

\newcommand{\ndims}{{d}} % n-dimensional
\newcommand{\domain}{\Omega} % domain
\newcommand{\imagedomain}{\Omega_{\rm scan}} % domain
\newcommand{\maxvoxel}{m_{\rm vox}} % number of voxels
\newcommand{\voxelsize}{\Delta} % voxel size
\newcommand{\basismeshsize}{h} % basis mesh size
\newcommand{\scanmesh}{\mathcal{T}_{\rm scan}^h}
\newcommand{\voxelmesh}{\mathcal{T}_{\rm scan}^\Delta}% basis mesh
\newcommand{\ambientdomain}{\Omega_{\rm scan}}
\newcommand{\basisorder}{k} % spline order
\newcommand{\regularity}{\alpha} % spline regularity
\newcommand{\order}{k} % spline order
\newcommand{\splinespace}{\mathcal{S}^{\order}_{\regularity}} % spline space
\newcommand{\basisfunc}{N} % basis function
\newcommand{\gradn}{\partial_n}
\renewcommand{\sup}{\operatornamewithlimits{sup\vphantom{p}}}

\newcommand{\boundary}{\partial \domain}
\newcommand{\dirichletboundary}{\boundary_{D}}
\newcommand{\neumannboundary}{\boundary_{N}}
\newcommand{\meshsize}{h}
\newcommand{\mesh}{\mathcal{T}^{\meshsize}}
\newcommand{\cutmesh}{\mesh_{\domain}}      % mesh with cut elements
\newcommand{\ambientmesh}{\mesh_{\ambientdomain}} % ambient mesh
\newcommand{\boundarymesh}{\mesh_{\boundary}}

\newcommand{\element}{K}
\newcommand{\face}{F}
\newcommand{\edge}{E}
\newcommand{\skeleton}{\mathcal{F}_{\rm skeleton}}
\newcommand{\ghost}{\mathcal{F}_{\rm ghost}}

\renewcommand{\u}{u}		% trail function space
\newcommand{\uh}{{\u^{\meshsize}}}		% trail function space
\newcommand{\bodystokes}{\boldsymbol{f}}
\newcommand{\neumanndatastokes}{\boldsymbol{t}}
\newcommand{\dirichletdatastokes}{\boldsymbol{g}}

\newcommand{\rintu}{\boldsymbol{r}_{{\rm int},\uu}^{\meshsize}}
\newcommand{\rintp}{r_{{\rm int},\p}^{\meshsize}}

\newcommand{\rjumpu}{\boldsymbol{r}_{\rm jump}^{\meshsize}}

\newcommand{\rghostu}{\boldsymbol{r}_{\rm ghost}^{\meshsize}}

\newcommand{\rnitscheu}{\boldsymbol{r}_{\rm nitsche}^{\meshsize}}
\newcommand{\rneumannu}{{\boldsymbol{r}_{\rm neumann}^{\meshsize}}}
\newcommand{\rskeleton}{r_{\rm skeleton}^{\meshsize}}

\newcommand{\uu}{{\boldsymbol{u}}}
\newcommand{\uuh}{{{\uu}^{\meshsize}}}
\renewcommand{\vv}{{\boldsymbol{v}}}
\newcommand{\vvh}{{{\vv}^{\meshsize}}}
\newcommand{\p}{p}
\newcommand{\ph}{{{\p}^{\meshsize}}}
\newcommand{\q}{q}
\newcommand{\qh}{{{\q}^{\meshsize}}}
\newcommand{\nn}{\boldsymbol{n}}			% normal vector
\newcommand{\ltrivert}{\left\vert\kern-0.25ex\left\vert\kern-0.25ex\left\vert}
\newcommand{\rtrivert}{\right\vert\kern-0.25ex\right\vert\kern-0.25ex\right\vert}

\usepackage[style=numeric,sorting=none,natbib=true,firstinits,citestyle=numeric-comp,doi=false,isbn=false,url=false,eprint=false,date=year,maxbibnames=2]{biblatex}
\renewbibmacro{in:}{%
  \ifentrytype{article}{}{\printtext{\bibstring{in}\intitlepunct}}}
\AtEveryBibitem{\clearfield{number}}
\begin{document}

\title*{Scan-based immersed isogeometric flow analysis}
\author{Clemens V.\ Verhoosel, E.\ Harald van Brummelen, Sai C.\ Divi and Frits de Prenter}
\authorrunning{Verhoosel \emph{et al.}}
\institute{Clemens V.\ Verhoosel \at Eindhoven University of Technology, PO Box 513, 5600 MB Eindhoven, The Netherlands \email{c.v.verhoosel@tue.nl}
\and E. Harald van Brummelen \at Eindhoven University of Technology, PO Box 513, 5600 MB Eindhoven, The Netherlands \email{e.h.v.brummelen@tue.nl}
\and Sai C.\ Divi \at Eindhoven University of Technology, PO Box 513, 5600 MB Eindhoven, The Netherlands \email{s.c.divi@tue.nl}
\and Frits de Prenter \at Delft University of Technology, PO Box 5, 2600 AA Delft, The Netherlands \email{f.deprenter@tudelft.nl}}

\maketitle

\abstract{This chapter reviews the work conducted by our team on scan-based immersed isogeometric analysis for flow problems.  To leverage the advantageous properties of isogeometric analysis on complex scan-based domains, various innovations have been made: \emph{(i)} A spline-based segmentation strategy has been developed to extract a geometry suitable for immersed analysis directly from scan data; \emph{(ii)} A stabilized equal-order velocity-pressure formulation for the Stokes problem has been proposed to attain stable results on immersed domains; \emph{(iii)} An adaptive integration quadrature procedure has been developed to improve computational efficiency; \emph{(iv)} A mesh refinement strategy has been developed to capture small features at \emph{a priori} unknown locations, without drastically increasing the computational cost of the scan-based analysis workflow. We review the key ideas behind each of these innovations, and illustrate these using a selection of simulation results from our work. A patient-specific scan-based analysis case is reproduced to illustrate how these innovations enable the simulation of flow problems on complex scan data.}

\footnotetext{This is a preprint of the following chapter: Clemens V. Verhoosel, E. Harald van Brummelen, Sai C. Divi and Frits de Prenter, Scan-based isogeometric flow analysis, published in Frontiers in Computational Fluid-Structure Interaction and Flow Simulation: Research from Lead Investigators under 40 - 2022, edited by Tayfun Tezduyar, 2022, Springer Nature Switzerland AG.}

\section{Introduction}
\label{sec:introduction}

\begin{figure}
  \centering
  \includegraphics[width=\textwidth]{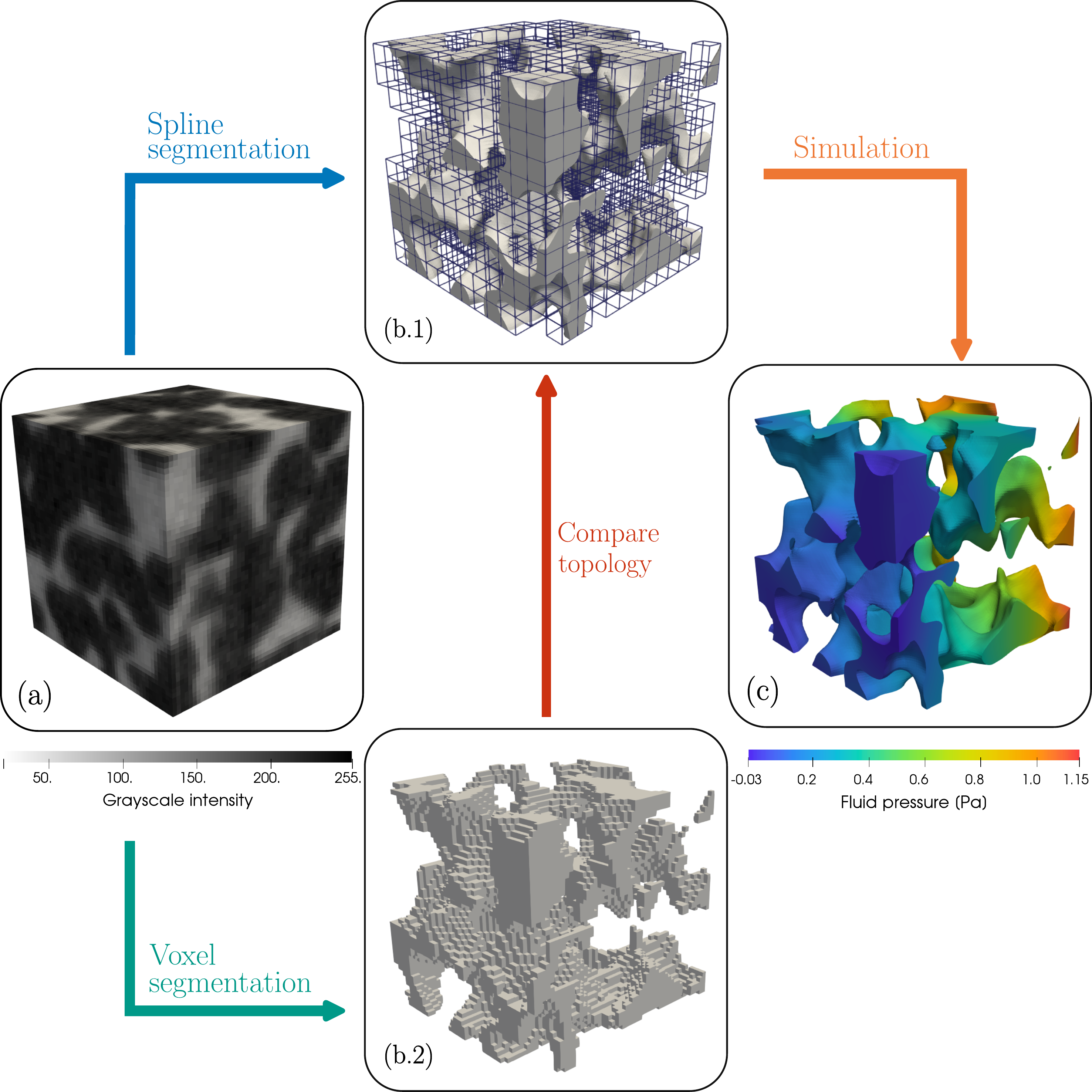}
  \caption{Illustration of the developed scan-based simulation workflow, considering a sintered glass specimen \cite{hoang_skeleton-stabilized_2019} as a typical example. The grayscale scan data is shown in panel (a). A smooth reconstruction of the geometry (representing the void space) based on a spline segmentation is shown in panel (b.1), along with the computational mesh in which the geometry is submersed. The directly segmented voxel image, which is used to assess the topological correctness of the spline segmentation, is shown in panel (b.2). A typical immersed isogeometric analysis result is shown in panel (c).}
  \label{fig:illustration}
\end{figure}

The rapid developments in the field of scientific computing have opened the doors to performing computational analyses on data obtained using advanced scanning technologies (\emph{e.g.}, tomography or photogrammetry). Such analyses are of particular interest in applications pertaining to non-engineered systems, which are common in, for example, biomechanics, geomechanics and material science. For scan-based simulations, the data sets from which the geometric models are constructed are typically very large, and the obtained models can be very complex in terms of both geometry and topology (see Fig.~\ref{fig:illustration}). In the context of standard finite element analyses (FEA), scan-based simulations require image segmentation and meshing techniques to produce high-quality analysis-suitable meshes that fit to the boundaries of the domain of interest. The construction of a FEA-suitable computational domain can be an error-prone and laborious process, involving manual geometry clean-up and mesh repairing and optimization operations. Such operations can account for the majority of the total computational analysis time and form a bottleneck in the automation of scan-based simulation workflows \cite{zhang_challenges_2013}.

The challenges associated with the simulation workflow for complex problems sparked the development of the isogeometric analysis (IGA) paradigm by Hughes and co-workers in 2005 \cite{hughes_isogeometric_2005}. The pivotal idea of IGA is to directly employ the geometry interpolation functions used in computer-aided design (\emph{e.g.}, B-splines and NURBS \cite{rogers_introduction_2001}) for the discretization of boundary value problems, thereby circumventing the problems associated with meshing. Besides the advantage of avoiding the meshing procedure and eliminating mesh-approximation errors, the use of higher-order continuous splines for the approximation of the solution has been demonstrated to yield accurate results using relatively few degrees of freedom for many (smooth) problems (see Ref.~\cite{cottrell_isogeometric_2009} for an overview). While isogeometric analysis has been successfully applied to complex three-dimensional problems based on (multi-patch) CAD objects (see, \emph{e.g.}, Refs.~\cite{cottrell_isogeometric_2006,zhang_patient-specific_2007,ruess_weak_2014,yu_anatomically_2020,hughes_chapter_2021,bucelli_multipatch_2021}), its application to scan-based simulations is hindered by the absence of analysis-suitable spline-based geometry models. Although spline preprocessors have been developed over the years for a range of applications \cite{zhang_solid_2012,hsu_interactive_2015,urick_review_2019}, the robust generation of analysis-suitable boundary-fitting volumetric splines for scan-based analyses is beyond the scope of the current tools on account of the geometrical and topological complexity typically inherent to scan data.

To still leverage the advantageous approximation properties of splines in scan-based simulations, IGA is often used in combination with immersed methods. In immersed methods, a non-boundary-fitting mesh is considered, in which the computational domain is submersed. Since the immersed domain does not align with the computational grid, some of the elements in the grid are cut by the immersed boundary and require a special treatment. The immersed approach has been considered in the finite element setting in the context of the Finite Cell Method (FCM) \cite{parvizian_finite_2007,duster_finite_2008,schillinger_finite_2015} and CutFEM \cite{burman_ghost_2010,burman_fictitious_2012,burman_cutfem_2015}, amongst others. The immersed concept has also been used in combination with IGA \cite{rank_geometric_2012,schillinger_isogeometric_2012,ruess_weakly_2013}, a strategy which is sometimes referred to as immersogeometric analysis \cite{kamensky_immersogeometric_2015,hsu_direct_2016}. The versatility of immersed isogeometric analysis techniques with respect to the geometry representation -- in the sense that the analysis procedure is not strongly affected by the complexity of the physical domain -- makes it particularly attractive in the scan-based analysis setting. Applications can nowadays be found in, for example, the modeling of trabecular bone \cite{verhoosel_image-based_2015,ruess_finite_2012,de_prenter_multigrid_2020}, porous media \cite{hoang_skeleton-stabilized_2019}, coated metal foams \cite{duster_numerical_2012}, metal castings \cite{jomo_robust_2019} and additive manufacturing \cite{carraturo_modeling_2020}.

Over the past decade our team has developed an analysis workflow using the  immersed isogeometric analysis paradigm. This workflow is illustrated in Fig.~\ref{fig:illustration}. In this chapter we review the key research contributions that made this analysis workflow applicable to scan-based flow simulations:
\begin{itemize}
    \item A \textbf{spline-based geometry segmentation} technique was proposed in Verhoosel \emph{et al.} \cite{verhoosel_image-based_2015}, with further improvements being made by Divi~\emph{et al.} \cite{divi_error-estimate-based_2020,divi_topology-preserving_2022}. The pivotal idea of the developed segmentation strategy is that the original scan data, which is usually non-smooth (\emph{i.e.}, a voxel representation), is smoothed using a spline approximation. A segmentation procedure able to provide an accurate explicit parametrization of the smoothed geometry then provides a geometric description of the scan object which is suitable for immersed isogeometric analysis.
    \item A \textbf{stabilized immersed isogeometric analysis} formulation for flow problems was proposed by Hoang \emph{et al.} \cite{hoang_skeleton-stabilized_2019}. The key idea of the proposed formulation is to use face-based stabilization techniques to make immersed simulations robust with respect to (unfavorably) cut elements, preventing the occurrence of oscillations in the velocity and pressure approximations. The stabilization terms also enable the consideration of equal-order discretizations of the velocity and pressure fields, which would otherwise cause inf-sup stability problems even in boundary-fitting finite elements \cite{hoang_mixed_2017}.
    \item An \textbf{adaptive integration procedure} was developed by Divi \emph{et al.} \cite{divi_error-estimate-based_2020} to reduce the computational cost involved in the evaluation of integrals over  cut elements, thereby improving the computational efficiency of the immersed analysis workflow. Based on Strang's lemma \cite{strang_analysis_2008}, an estimator for the integration error is derived, which is then used to optimally distribute integration quadrature points over cut elements.
    \item An \textbf{error-estimation-based adaptive refinement} procedure has been developed by Divi \emph{et al.} \cite{divi_residual-based_2022} to capture small features without drastically increasing the computational cost of the scan-based workflow. Residual-based error estimators are constructed to perform local basis function refinements to increase the resolution of the spline basis in regions where this is particularly beneficial from an accuracy point of view, without prior knowledge of the locations of these regions.
\end{itemize}
Our scan-based immersed isogeometric analysis workflow has been applied to a range of real world data problems, mainly in the context of $\mu$CT-scans. In this chapter we illustrate the capabilities of our workflow in the context of patient-specific arterial flow problems. The analysis of porous medium flows as presented in Ref.~\cite{hoang_skeleton-stabilized_2019}, and illustrated in Fig.~\ref{fig:illustration}, forms another prominent application of our method. 

This chapter is organized as follows. The essential innovations regarding each of the research contributions listed above are reviewed in Sections~\ref{sec:segmentation}--\ref{sec:adaptivesplines}. A typical application of the developed workflow will then be discussed in Section~\ref{sec:results}. We will conclude this chapter in Section~\ref{sec:conclusions} with an assessment of our scan-based analysis workflow, discussing its capabilities and current limitations.

\section{Spline-based geometry segmentation}
\label{sec:segmentation}

In this section we review the spline-based image segmentation procedure that we have developed in the context of scan-based immersed isogeometric analysis \cite{verhoosel_image-based_2015,divi_error-estimate-based_2020,divi_topology-preserving_2022}. In Section~\ref{sec:levelset} we first discuss the spline-based level set construction to smoothen scan data. In Section~\ref{sec:tessellation} we review the algorithms used to construct an explicit parametrization of the scan domain. An example is finally shown in Section~\ref{sec:topology}, illustrating the effectivity of the topology-preservation procedure developed in Ref.~\cite{divi_topology-preserving_2022}.

\begin{figure}
    \centering
    \includegraphics[width=\textwidth]{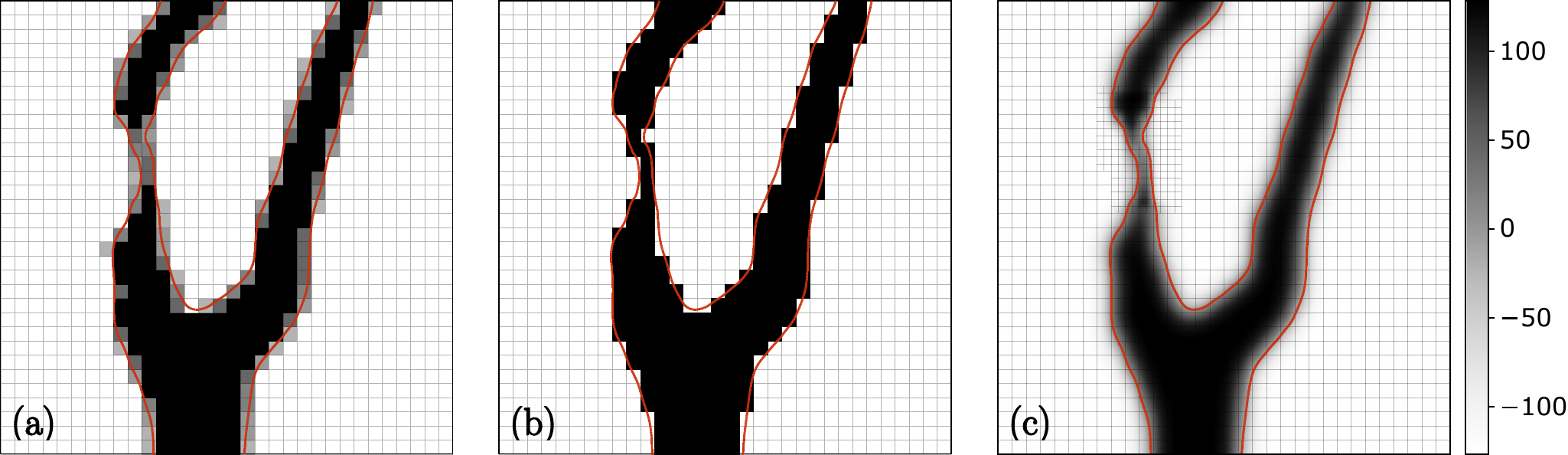}
    \caption{A two-dimensional illustration of the spline-based segmentation procedure. (a) Grayscale data on a $32 \times 32$ voxel mesh, $\voxelmesh$, as in \eqref{eq:grayscale}. (b) The voxel segmentation obtained by thresholding the grayscale data as in \eqref{eq:voxeldomain}. (c) The smooth level set function obtained using a THB-spline basis constructed on the voxel mesh with local refinements as in \eqref{eq:bsplinefunc}. The exact immersed boundary, which is in principle unknown, is shown in red for reference.}
    \label{fig:voxels}
\end{figure}

\subsection{B-spline smoothing of the scan data}\label{sec:levelset}
The spline-based level set construction is illustrated in Fig.~\ref{fig:voxels}. We consider a $\ndims$-dimensional scan domain, $\imagedomain= [0,L_{1}] \times\ldots \times [0,L_{\ndims}]$ with volume $V_{\rm scan}= \prod_{i=1}^{n_d} L_i$, which is partitioned by a set of $\maxvoxel$ voxels, as illustrated in Fig.~\ref{fig:voxels}a. We denote the voxel mesh by $\voxelmesh$, with $\boldsymbol{\voxelsize}$ the voxel size in each direction. The grayscale intensity function is then defined as
\begin{equation}
g: \voxelmesh \rightarrow \mathscr{G},
\label{eq:grayscale}
\end{equation}
with $\mathscr{G}$ the range of the grayscale data (\emph{e.g.}, from 0 to 255 for 8 bit unsigned integers). An approximation of the object $\domain$ can be obtained by thresholding the grayscale data,
\begin{equation}
  \domain\approx \{ \boldsymbol{x} \in \imagedomain | g(\boldsymbol{x}) > g_{\rm crit} \}\subset \imagedomain,
  \label{eq:voxeldomain}
\end{equation}
where $g_{\rm crit}$ is the threshold value. As a consequence of the piecewise definition of the grayscale data in equation \eqref{eq:grayscale}, the boundary of the segmented object is non-smooth when the grayscale data is segmented directly. In the context of analysis, the non-smoothness of the boundary can be problematic, as irregularities in the surface may lead to non-physical features in the solution to the problem.

The spline-based segmentation procedure developed in Refs.~\cite{verhoosel_image-based_2015,divi_error-estimate-based_2020,divi_topology-preserving_2022} enables the construction of a smooth boundary approximation based on voxel data. The key idea of this spline-based segmentation technique is to smoothen the grayscale function \eqref{eq:grayscale} by convoluting it using an $n$-dimensional spline basis, $\{ \basisfunc_{i,\basisorder}(\scanmesh) \}_{i=1}^n$, defined over a mesh, $\scanmesh$, with element size, $\boldsymbol{h}=(h_1,\ldots,h_d)$ (note that the mesh size can differ from the voxel size). The order $\basisorder$ of the spline basis functions is assumed to be constant and isotropic. We consider THB-splines \cite{giannelli_thb-splines_2012} for the construction of locally refined spaces. By considering full-regularity ($C^{\basisorder-1}$-continuous) splines of degree $\basisorder>1$, a smooth level set approximation of \eqref{eq:grayscale} is obtained by the convolution operation
\begin{align}
f(\boldsymbol{x}) &= \sum \limits_{i=1}^{n} \basisfunc_{i,\basisorder} (\boldsymbol{x}) a_{i}, &  a_{i} &= \cfrac{\int_{\imagedomain} \basisfunc_{i,\basisorder} (\boldsymbol{x}) g(\boldsymbol{x}) {\rm d}\boldsymbol{x}}{\int_{\imagedomain} \basisfunc_{i,\basisorder} (\boldsymbol{x}) {\rm d}\boldsymbol{x}}  , \label{eq:bsplinefunc}
\end{align}
where $\{ a_i \}_{i=1}^n$ are the coefficients of the discrete level set function. The smoothed domain then follows by thresholding of this level set function:
\begin{equation}
  \domain\approx \{ \boldsymbol{x} \in \imagedomain | f(\boldsymbol{x}) > f_{\rm crit} \}\subset \imagedomain.
  \label{eq:implicitdomain}
\end{equation}

The spline level set function corresponding to the voxel data in Fig.~\ref{fig:voxels}a is illustrated in Fig.~\ref{fig:voxels}c for the case of a locally refined mesh $\scanmesh$ and second order ($\basisorder=2$) THB-splines. As can be seen, the object retrieved from the convoluted level set function more closely resembles the original geometry in Fig.~\ref{fig:voxels}a compared to the voxel segmentation in Fig.~\ref{fig:voxels}b. Also, as a consequence of the higher-order continuity of the spline basis, the boundaries of the domain are smooth, which is in closer agreement with reality.

The convolution operation \eqref{eq:bsplinefunc} is computationally efficient, resulting from the fact that it is not required to solve a linear system of equations (in contrast to a (global) $L^2$-projection) and the restricted support of the convolution kernel. Moreover, the convolution strategy has various properties that are advantageous in the context of scan-based immersed isogeometric analysis (see Refs.~\cite{verhoosel_image-based_2015,divi_topology-preserving_2022} for details):

\subsubsection*{Conservation of the gray scale intensity}
Under the condition that the spline basis, $\{ \basisfunc_{i,\basisorder}(\scanmesh) \}_{i=1}^n$, satisfies the partition of unity property (\emph{e.g.}, B-splines, THB-splines), the smooth level set approximation \eqref{eq:bsplinefunc} conserves the gray scale intensity of the original data in the sense that
\begin{align}
\frac{1}{V_{\rm scan}} \int \limits_{\Omega_{\rm scan}} f\,{\rm d}V =  \frac{1}{V_{\rm scan}} \int \limits_{\Omega_{\rm scan}} g\,{\rm d}V = \frac{1}{m_{\rm vox}} \sum_{\element \in \voxelmesh} g(\element).
 \label{eq:intensity}
\end{align}
This property ensures that there is a direct relation between the threshold value, $f_{\rm crit}$, for the smooth level set reconstruction \eqref{eq:implicitdomain} and that of the original data, $g_{\rm crit}$, in equation \eqref{eq:voxeldomain}.

\subsubsection*{Local boundedness by the original data}
On every voxel $\element \in \voxelmesh$, the level set function \eqref{eq:bsplinefunc} is bound by the extrema of the voxel function over the support extension \cite{bazilevs_isogeometric_2006}, $\widetilde{\element}$, \emph{i.e.},
\begin{align}
 &\min_{\boldsymbol{x} \in  \widetilde{\element}} g(\boldsymbol{x})  \leq f(\boldsymbol{x}) \leq \max_{\boldsymbol{x} \in \widetilde{\element}} g(\boldsymbol{x}) & &\forall \boldsymbol{x} \in \element.
 \label{eq:bounds}
\end{align}
These bounds preclude overshoots and undershoots, which indicates that no spurious oscillations are created by the smoothing procedure (contrasting the case of an $L^2$-projection).

\subsubsection*{Approximate Gaussian blurring}
The spline-based convolution operation \eqref{eq:bsplinefunc} can be written as an integral transform
\begin{align}
f(\boldsymbol{x}) &= \int_{\imagedomain}\mathcal{K}(\boldsymbol{x},\boldsymbol{y}) g(\boldsymbol{y}) \:{\rm d} \boldsymbol{y},  &  \mathcal{K}(\boldsymbol{x}, \boldsymbol{y}) &= \sum \limits_{i=1}^{n} \cfrac{\basisfunc_{i,\basisorder}(\boldsymbol{x}) \basisfunc_{i,\basisorder}(\boldsymbol{y})}{\int_{\imagedomain} \basisfunc_{i,\basisorder}(\boldsymbol{z})\:{\rm d}\boldsymbol{z} }, \label{eq:convolution}
\end{align}
where $\mathcal{K}(\boldsymbol{x},\boldsymbol{y})$ is the kernel of the transformation.

\begin{figure}
	\centering
	\begin{subfigure}[t]{0.48\textwidth}
		\centering
		\includegraphics[width=\textwidth]{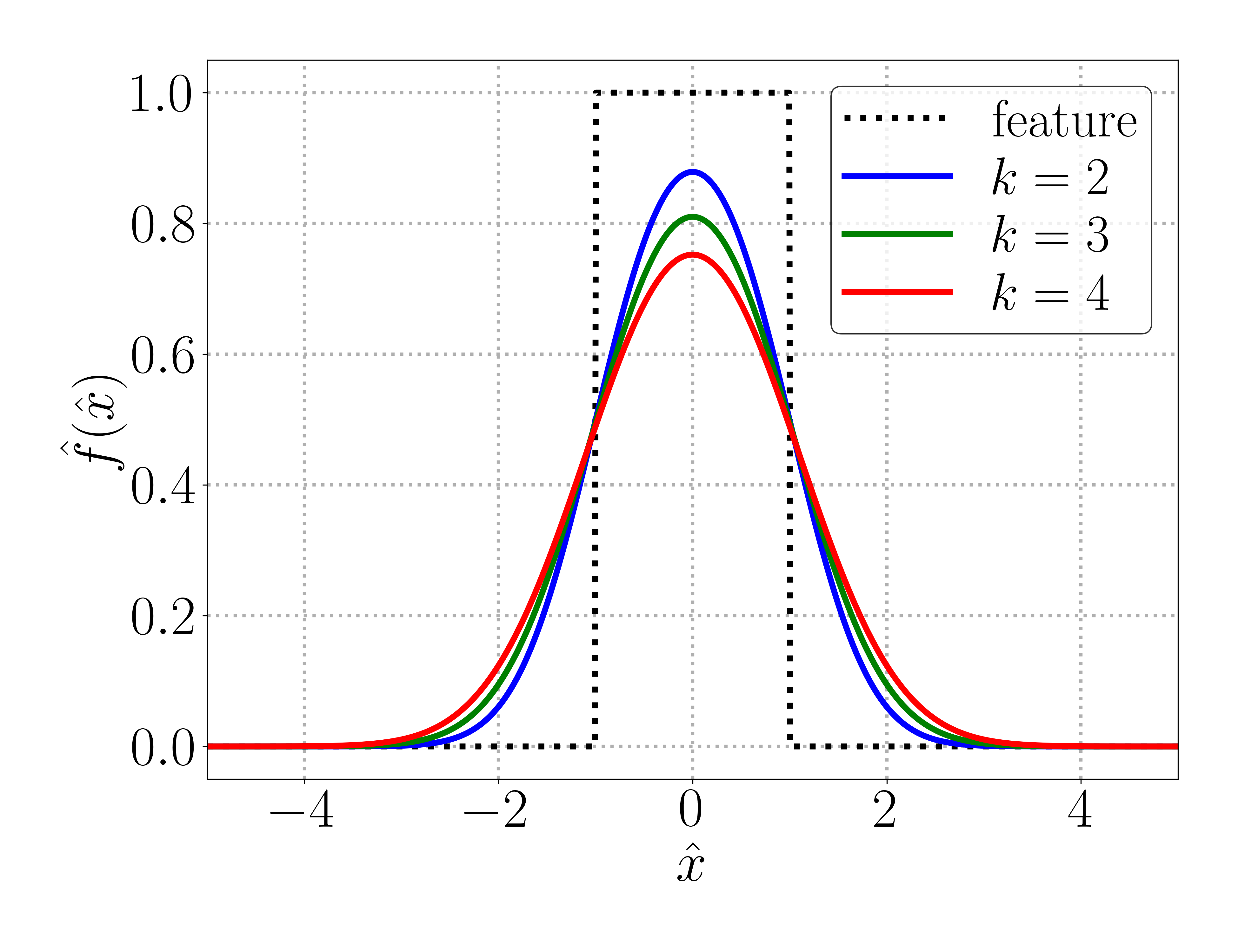}
		\caption{$\hat{\ell}=\frac{\ell}{h}=2$}
		\label{fig:spaceloverh2}
	\end{subfigure}\\
	\begin{subfigure}[t]{0.48\textwidth}
		\centering
		\includegraphics[width=\textwidth]{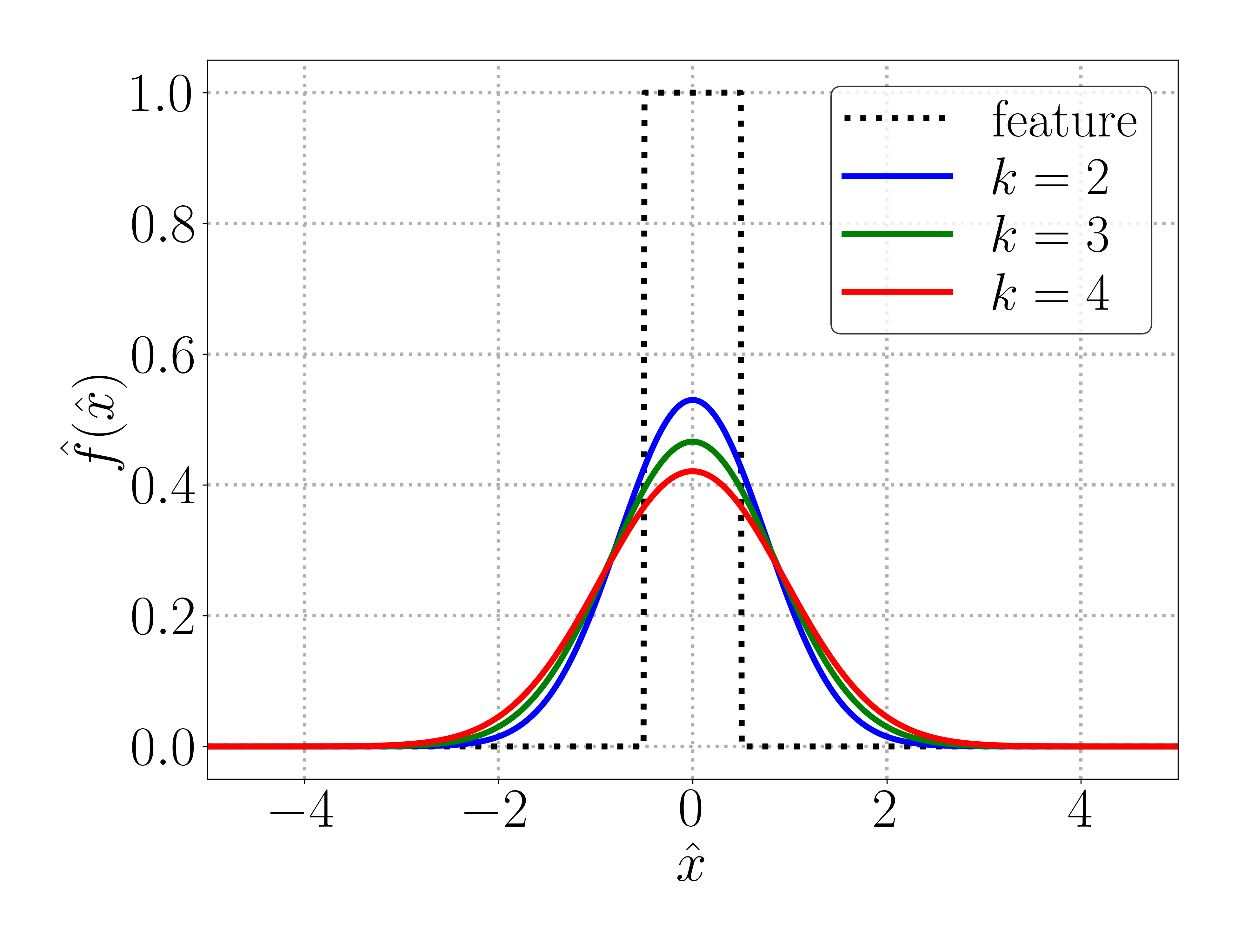}
		\caption{$\hat{\ell}=\frac{\ell}{h}=1$}
		\label{fig:spaceloverh1}
	\end{subfigure}\hfill%
	\begin{subfigure}[t]{0.48\textwidth}
		\centering
		\includegraphics[width=\textwidth]{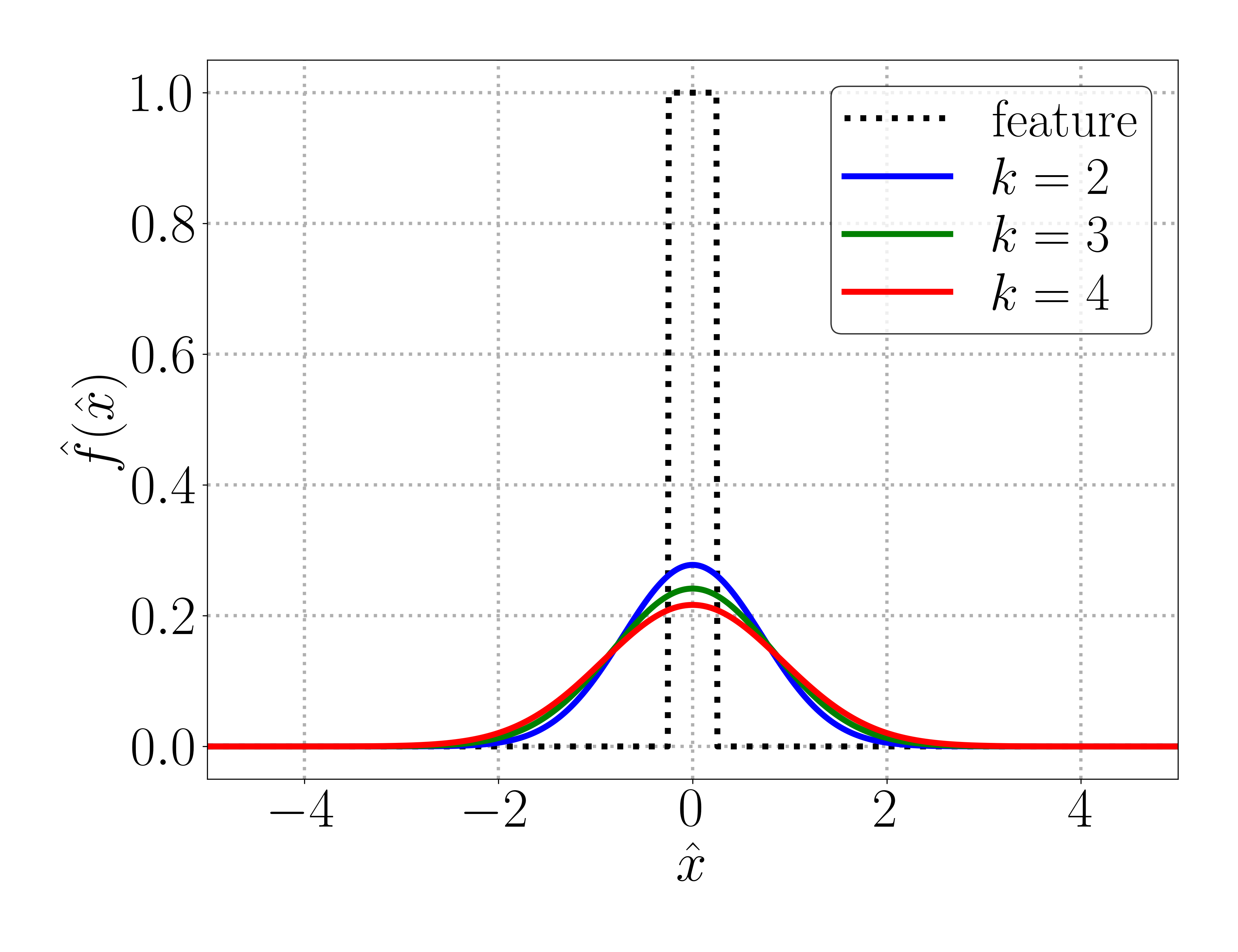}
		\caption{$\hat{\ell}=\frac{\ell}{h}=\frac{1}{2}$}
		\label{fig:spaceloverhhalf}
	\end{subfigure}%
	\caption{Smoothed level set approximation \eqref{eq:bsplinefunc} of a geometric feature in the spatial domain, with $\hat{x}=x/h$, for various feature-size-to-mesh ratios, $\hat{\ell}=\ell/h$, and B-spline degrees.}
	\label{fig:levelsetapproximation}
\end{figure}

The integral transform \eqref{eq:convolution} acts as an approximate Gaussian filter \cite{deng_adaptive_1993}. We illustrate this behavior for the case of one-dimensional voxel data, which is smoothed using a B-spline basis defined on a uniform mesh, $\scanmesh$, with mesh size $\basismeshsize$. Following the derivation in Ref.~\cite{verhoosel_image-based_2015} -- in which the essential step is to approximate the B-spline basis functions by rescaled Gaussians \cite{unser_asymptotic_1992} -- the integration kernel \eqref{eq:convolution} can be approximated by
\begin{equation}
\mathcal{K}(x,y) \approx \kappa( x-y ) = \frac{1}{\sigma \sqrt{2\pi}} \exp{\left( {-\cfrac{(x-y)^2}{2 \sigma^2}} \right)}, \label{eq:kernel}
\end{equation}
where the width of the smoothing kernel is given by $\sigma = \basismeshsize \sqrt{ \frac{\basisorder+1}{6} }$. Next, we consider an object of size $\ell$, represented by the grayscale function
\begin{align}
 g(x) &= \begin{cases} 
   1 & | x |  < \ell/2\\
  0 & \mbox{otherwise}
\end{cases}.
\end{align}
This object and the corresponding approximate level set function \eqref{eq:bsplinefunc} are illustrated in Fig.~\ref{fig:levelsetapproximation} for various feature-size-to-mesh ratios, $\hat{\ell}=\ell/\basismeshsize=2,1,\frac{1}{2}$, and B-spline degrees, $\basisorder=2,3,4$. Following the (Fourier) analysis in Ref.~\cite{divi_topology-preserving_2022}, the value of the smoothed level set function at $x=0$ follows as
\begin{equation}
 \hat{f}_1(0) = \hat{\ell} \sqrt{\frac{3}{\pi(\basisorder+1)} }  ~ \exp{\left( - \frac{3 \hat{\ell}^2}{16(\basisorder+1)}   \right)},
\label{eq:approxmax}
\end{equation}
which conveys that the maximum value of the smoothed level set depends linearly on the relative feature size $\hat{\ell}$ (for sufficiently small $\hat{\ell}$), and decreases with increasing B-spline order.

Fig.~\ref{fig:spaceloverh2} shows the case for which the considered feature is twice as large as the mesh size, \emph{i.e.}, $\hat{\ell}=2$, illustrating that the sharp boundaries of the original grayscale function are significantly smoothed. The decrease in the maximum level set value as given by equation \eqref{eq:approxmax} is observed. When the level set function is segmented by a threshold of $g_{\rm crit}=0.5$, a geometric feature that closely resembles the original one is recovered. Fig.~\ref{fig:spaceloverh1}-\subref{fig:spaceloverhhalf} illustrate cases where the feature length, $\ell$, is not significantly larger than the mesh size, $\basismeshsize$. For the case where the feature size is equal to the size of the mesh, the maximum of the level set function drops significantly compared to the case of $\hat{\ell}=2$. When considering second-order B-splines, the maximum is still marginally above $g_{\rm crit}=0.5$. Although the recovered feature is considerably smaller than the original one, it is still detected in the segmentation procedure. When increasing the B-spline order, the maximum value of the level set drops below the segmentation threshold, however, indicating that the feature will no longer be detected. When decreasing the feature length further, as illustrated in Fig.~\ref{fig:spaceloverhhalf}, the feature would be lost when segmentation is performed with $g_{\rm crit}=0.5$, regardless of the order of the spline basis. 

The implications of this smoothing behavior of the convolution operation \eqref{eq:bsplinefunc} in the context of the spline-based segmentation will be discussed in Section~\ref{sec:topology}.

\subsection{Octree-based tessellation procedure} \label{sec:tessellation}
\begin{figure}
    \centering
    \includegraphics[width=\textwidth]{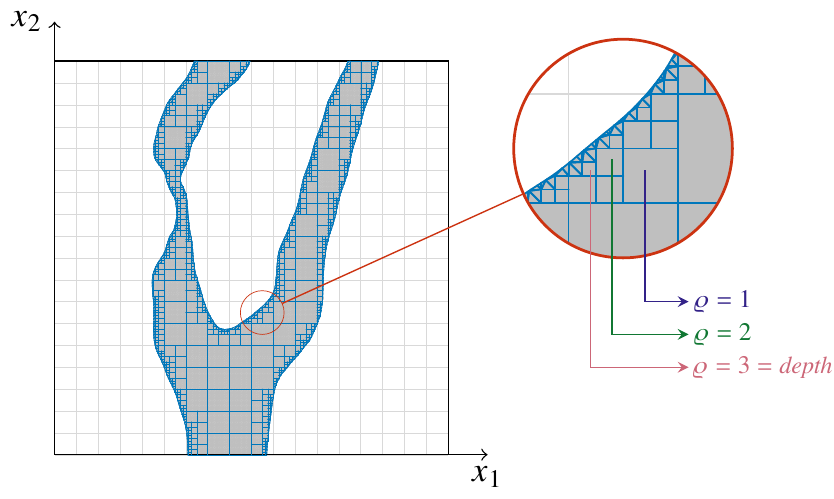}
    \caption{Illustration of the octree-based tessellation procedure to acquire an explicit parametrization of the immersed computational domain. The number of element subdivisions is indicated by $\varrho$, which is equal to the octree \textit{depth} at the lowest level of bisectioning.}
    \label{fig:octree}
\end{figure}

Our scan-based isogeometric analysis approach requires the construction of an explicit parametrization of the implicit level set domain \eqref{eq:implicitdomain}. In this section we outline the segmentation procedure that we use to obtain an explicit parametrization of the domain and its (immersed) boundaries. This procedure, which is based on the octree subdivision approach introduced in the context of the Finite Cell Method in Ref.~\cite{duster_finite_2008}, is illustrated in Fig.~\ref{fig:octree}. In Section~\ref{sec:octree} we first discuss the employed octree procedure, after which the tessellation procedure used at the lowest level of subdivision is detailed in Section~\ref{sec:midpoint}. Without loss of generality, in the remainder we will assume that the level set function $f$ is shifted such that $f_{\rm crit}=0$.

\subsubsection{Octree subdivision} \label{sec:octree}
In our analysis framework we consider a regular mesh $\scanmesh$ that conforms to the scan domain $\imagedomain$. Each element in this mesh is a hyperrectangle (\emph{i.e.}, a line in one dimension, a rectangle in two dimensions, and a hexahedron in three dimensions) with size ${\boldsymbol{h}}=(h_1,\ldots, h_{\ndims})$ in each direction. Elements, $K$, that are intersected by the immersed boundary are trimmed, resulting in a partitioning, $\mathcal{P}_K$, of a cut element (see Fig.~\ref{fig:octree}).

\begin{algorithm}
\caption{Function that trims an element based on an evaluate level set function}\label{alg:trim}
    \hspace*{\algorithmicindent} \textbf{Input:} array of level set \textit{values}, octree \textit{depth}, \textit{dimension} of the element to be trimmed\\
    \hspace*{\algorithmicindent} \textbf{Output:} element of type $\textit{Hyperrectangle} \cup  \textit{Void} \cup \textit{WithChildren} \cup \textit{Mosaic}$
\begin{algorithmic}[1]
\Function{trim\_element}{\textit{values}, \textit{depth}, \textit{dimension}}\label{alg:trim:function}

\If{$\textbf{all}~\textit{values} > 0$}\Comment{All level set values are positive}\label{alg:trim:allpositive}
    \State \textbf{return} \textit{Hyperrectangle}(\textit{dimension})
\ElsIf{$\textbf{all}~\textit{values} \leq 0$}\Comment{All level set values are non-positive}\label{alg:trim:allnonpositive}
    \State \textbf{return} \textit{Void}(\textit{dimension})
\Else\label{alg:trim:mixed}
    \Comment{Both positive and non-positive level set values}
    \If{$\textit{depth} > 0$}
        \Comment{Recursively trim the element when $\textit{depth}>0$}

        \State $\textit{withchildren\_element} \gets \Call{get\_withchildren\_element}{\textit{dimension}}$\label{alg:trim:getwithchildren}
        \For{\textit{child} \textbf{in} \textit{withchildren\_element}}
            \State $\textit{child\_values}\gets \Call{get\_child\_values}{\textit{child}, \textit{values}}$\label{alg:trim:childvalues}
            \State \textit{child} $\gets$ \Call{trim\_element}{\textit{child\_values}, $\textit{depth}-1$, \textit{dimension}}\label{alg:trim:recursive}
        \EndFor
        \State \textbf{return} $\textit{withchildren\_element}$
    \Else\label{alg:trim:terminate}
        \Comment{Terminate the recursion when $\textit{depth}=0$}
        \State \textbf{return}  $\Call{get\_mosaic\_element}{\textit{values}, \textit{dimension}}$\label{alg:trim:mosaic}
    \EndIf
\EndIf
\EndFunction
\end{algorithmic}
\end{algorithm}

We employ the octree-based trimming procedure outlined in Alg.~\ref{alg:trim}. This procedure follows a bottom-up approach \cite{varduhn_tetrahedral_2016} in which the level set function \eqref{eq:bsplinefunc} is sampled at the $2^{\textit{depth}}+1$ vertices of the octree in each direction for each element, where $\textit{depth}$ is the number of subdivision operations performed to detect the immersed boundary.

The trimming procedure takes the evaluated level set \textit{values}, the subdivision \textit{depth} and the \textit{dimension} of the element as input arguments. If all level set \textit{values} are positive (L\ref{alg:trim:allpositive}), the \textsc{trim\_element} function retains the complete element in the mesh. If all level set values are non-positive (L\ref{alg:trim:allnonpositive}), the element is discarded. When some of the level set \textit{values} are positive and some are non-positive (L\ref{alg:trim:mixed}), this implies that the element is intersected by the immersed boundary. In this case, the element is subdivided (bisected) into $2^\ndims$ children (L\ref{alg:trim:getwithchildren}). For each \textit{child} a recursive call to the \textsc{trim\_element} function is made (L\ref{alg:trim:recursive}).

The recursive subdivision routine is terminated at the lowest level of subdivision, \emph{i.e.}, at $\textit{depth}=0$ (L\ref{alg:trim:terminate}). Since our analysis approach requires the evaluation of functions on the immersed boundary, it is convenient to also obtain an explicit parametrization of this boundary. To obtain this parametrization, at the lowest level of subdivision we consider a tessellation procedure based on the level set values. The function \textsc{get\_mosaic\_element} that implements this procedure is discussed in Section~\ref{sec:midpoint}.

\subsubsection{Midpoint tessellation}\label{sec:midpoint}
\begin{figure}
	\centering
	\begin{subfigure}[b]{0.32\textwidth}
		\centering
		\includegraphics[width=\textwidth]{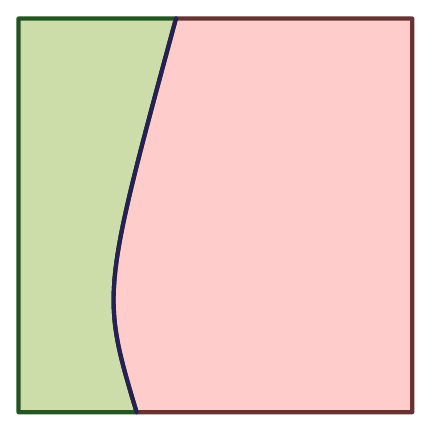}
		\caption{}
		\label{fig:cut_cell_2d}
	\end{subfigure}\hfill%
	\begin{subfigure}[b]{0.32\textwidth}
		\centering
		\includegraphics[width=\textwidth]{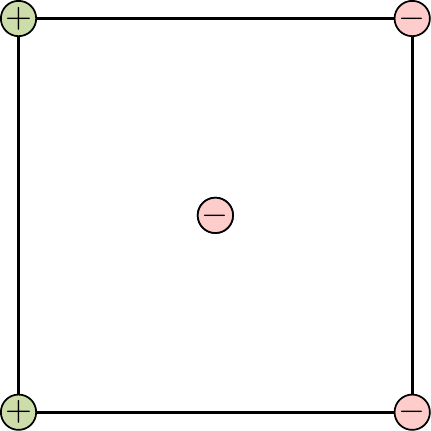}
		\caption{}
		\label{fig:levelset_vertex_2d}
	\end{subfigure}\hfill%
	\begin{subfigure}[b]{0.32\textwidth}
		\centering
		\includegraphics[width=\textwidth]{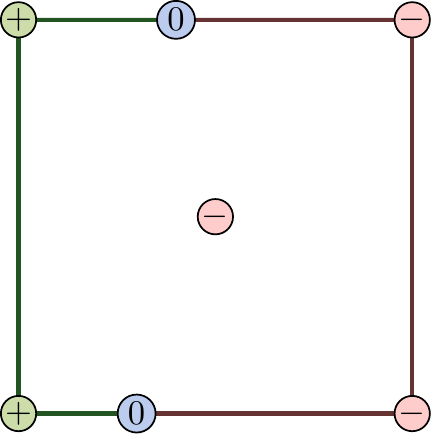}
		\caption{}
		\label{fig:trim_edges_2d}
	\end{subfigure}\\[12pt]
	\begin{subfigure}[b]{0.32\textwidth}
		\centering
		\includegraphics[width=\textwidth]{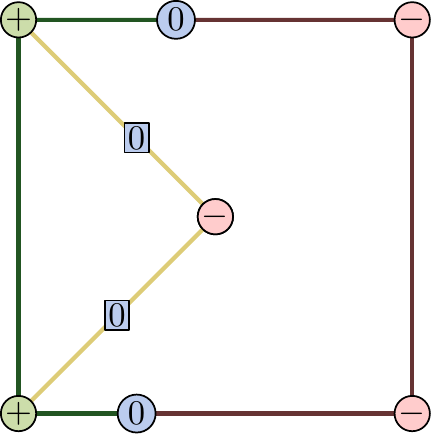}
		\caption{}
		\label{fig:levelset_diagonal_2d}
	\end{subfigure}\hfill%
	\begin{subfigure}[b]{0.32\textwidth}
		\centering
		\includegraphics[width=\textwidth]{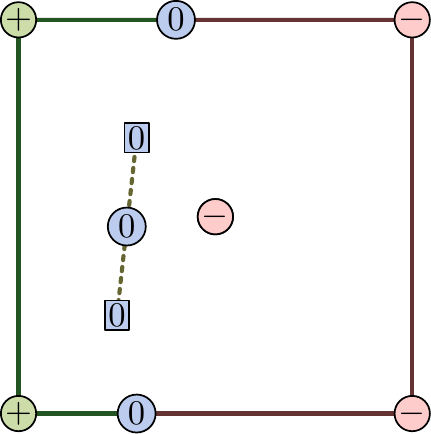}
		\caption{}
		\label{fig:diagonal_midpoint_2d}
	\end{subfigure}\hfill%
	\begin{subfigure}[b]{0.32\textwidth}
		\centering
		\includegraphics[width=\textwidth]{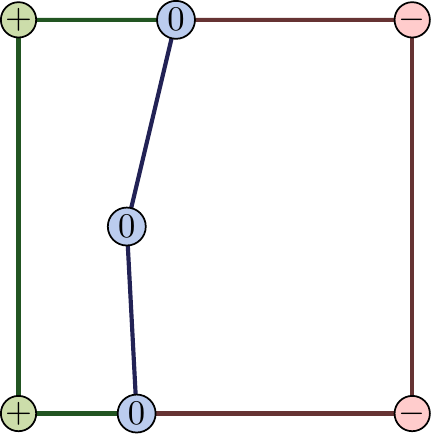}
		\caption{}
		\label{fig:levelset_midpoint_2d}
	\end{subfigure}\\[12pt]
	\begin{subfigure}[b]{0.32\textwidth}
		\centering
		\includegraphics[width=\textwidth]{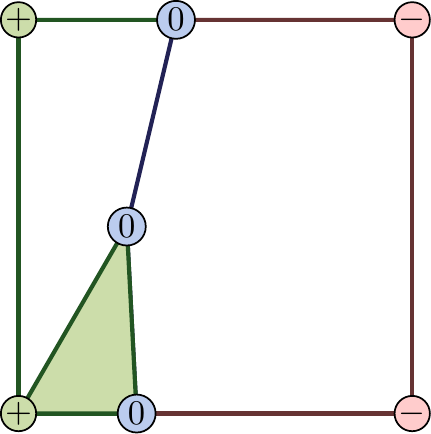}
		\caption{}
		\label{fig:tessellation1_2d}
	\end{subfigure}\hfill%
	\begin{subfigure}[b]{0.32\textwidth}
		\centering
		\includegraphics[width=\textwidth]{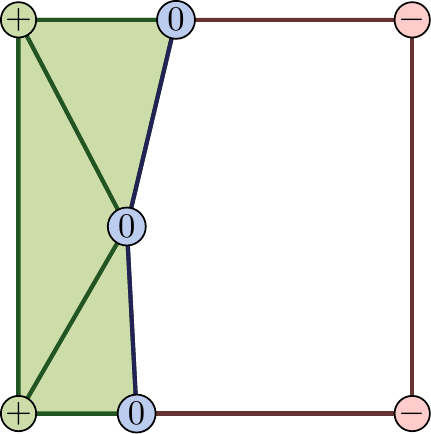}
		\caption{}
		\label{fig:tessellation2_2d}
	\end{subfigure}\hfill%
	\begin{subfigure}[b]{0.32\textwidth}
		\centering
		\includegraphics[width=\textwidth]{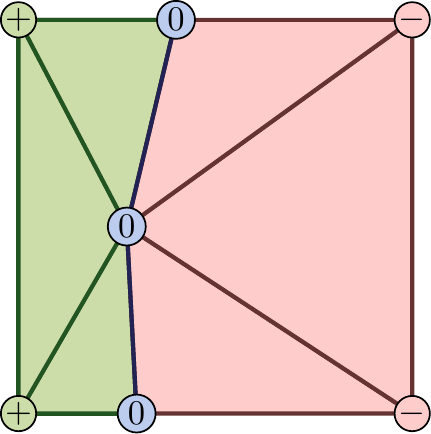}
		\caption{}
		\label{fig:tessellation3_2d}
	\end{subfigure}
	\caption{Schematic representation of the midpoint tessellation procedure for a two-dimensional case.}
	\label{fig:midpoint_2d}
\end{figure}

At the lowest level of subdivision of the octree procedure, we perform a tessellation based on the $2^\ndims$ level set values at that level. In the scope of our work, an important requirement for the tessellation procedure is that it is suitable for the consideration of interface problems. Practically, this means that if the procedure is applied to the negated level set function, a partitioning of the complementary part of the cell is obtained, with an immersed boundary that matches with that of the tessellation based on the original level set function (illustrated in Fig.~\ref{fig:tessellation3_2d}). Standard tessellation procedures, specifically Delaunay tessellation \cite{delaunay_sur_1934}, do not meet this requirement, as the resulting tessellation is always convex. When an immersed boundary tessellation is convex from one side, it is concave from the other side, meaning that it cannot be identically represented by the Delaunay tessellation. Another complication of Delaunay tessellation is its lack of uniqueness \cite{de_berg_computational_2008}, which in our applications results in non-matching interface tessellations.  

\begin{algorithm}
\caption{Function that constructs a mosaic element based on level set values}\label{alg:midpoint}
    \hspace*{\algorithmicindent} \textbf{Input:} array of level set \textit{values}, \textit{dimension} of the element to be constructed\\
    \hspace*{\algorithmicindent} \textbf{Output:} \textit{element} of type $\textit{Mosaic} \cup \textit{Hyperrectangle} \cup \textit{Void}$
\begin{algorithmic}[1]
\Function{get\_mosaic\_element}{$\textit{values}, \textit{dimension}$}
    \State \textit{hyperrectangle} $\gets$ \Call{get\_hyperrectangle}{\textit{dimension}}
    \If{$\textbf{all}~\textit{values} > 0$}\Comment{All level set values are positive}\label{alg:midpoint:allpositive}
        \State \textbf{return} \textit{hyperrectangle}
    \ElsIf{$\textbf{all}~\textit{values} \leq 0$}\label{alg:midpoint:allnonpositive}
        \Comment{All level set values are non-positive}
        \State \textbf{return} \textit{Void}
    \Else \label{alg:midpoint:mixed}
        \Comment{Both positive and non-positive level set values}
        \State \textit{boundaries} $\gets$ \Call{get\_boundaries}{\textit{hyperrectangle}}\label{alg:midpoint:startrecursion}
        \For{\textit{boundary} \textbf{in} \textit{boundaries}}
            \State $\textit{boundary\_values}\gets \Call{get\_boundary\_values}{\textit{boundary}, \textit{values}}$
            \State \textit{boundary} $\gets$ \Call{get\_mosaic\_element}{\textit{boundary\_values}, $\textit{dimension}-1$}\label{alg:midpoint:recursive}
        \EndFor\label{alg:midpoint:endrecursion}
        \State \textit{zero\_points} $\gets$ \Call{get\_zero\_points}{\textit{boundaries}}\label{alg:midpoint:zeropoints}
        \State \textit{midpoint} $\gets$ \Call{average}{\textit{zero\_points}}\label{alg:midpoint:midpoint}
        \State \textit{mosaic\_element} $\gets$ \Call{extrude\_to\_point}{\textit{boundaries}, \textit{midpoint}}\label{alg:midpoint:extrude}
        \State \textbf{return} \textit{mosaic\_element}
    \EndIf
\EndFunction
\end{algorithmic}
\end{algorithm}

\begin{figure}
	\centering
	\begin{subfigure}[b]{0.33\textwidth}
		\centering
		\includegraphics[width=\textwidth]{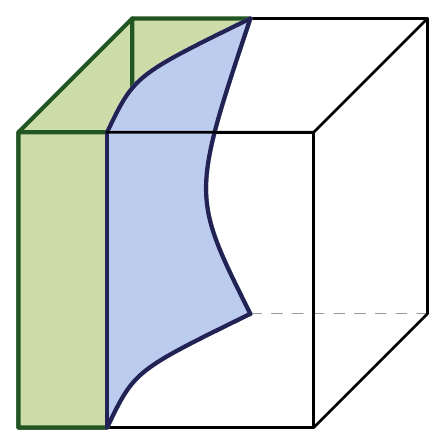}
		\caption{}
		\label{fig:cut_cell_3d}
	\end{subfigure}\hfill%
	\begin{subfigure}[b]{0.33\textwidth}
		\centering
		\includegraphics[width=\textwidth]{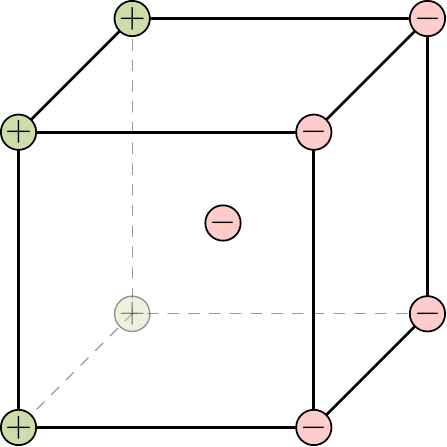}
		\caption{}
		\label{fig:levelset_vertex_3d}
	\end{subfigure}\hfill%
	\begin{subfigure}[b]{0.33\textwidth}
		\centering
		\includegraphics[width=\textwidth]{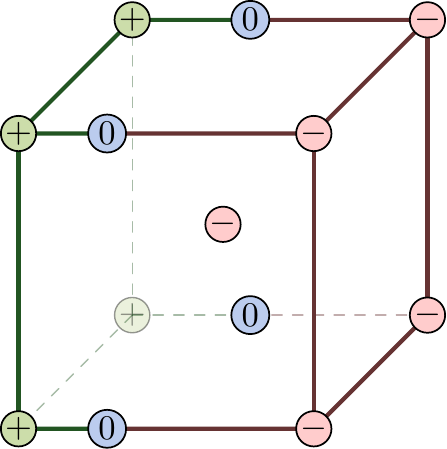}
		\caption{}
		\label{fig:trim_edges_3d}
	\end{subfigure}\\[12pt]
	\begin{subfigure}[b]{0.4\textwidth}
		\centering
		\includegraphics[width=\textwidth]{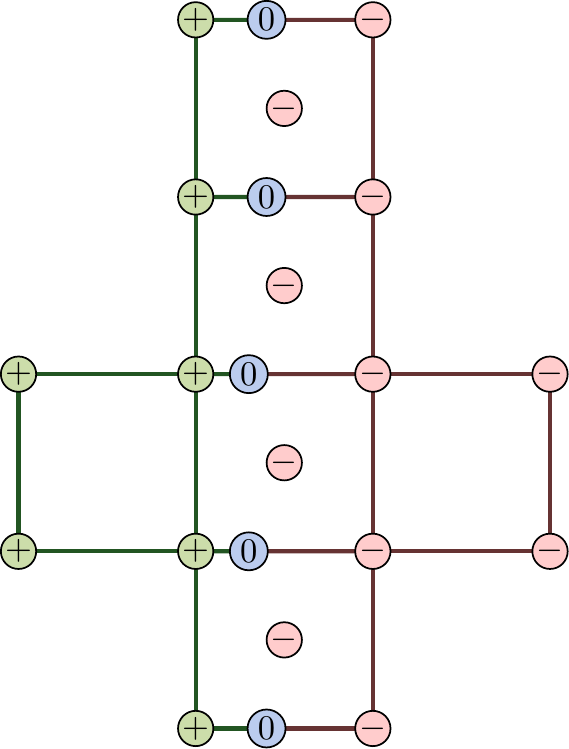}
		\caption{}
		\label{fig:diagonal_3d}
	\end{subfigure}\hfill%
		\begin{subfigure}[b]{0.4\textwidth}
		\centering
		\includegraphics[width=\textwidth]{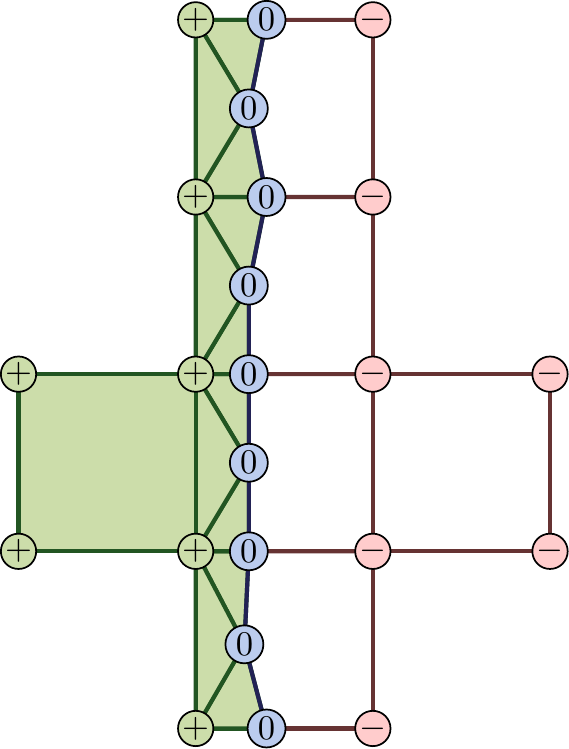}
		\caption{}
		\label{fig:diagonal_midpoint_3d}
	\end{subfigure}\\[12pt]
	\begin{subfigure}[b]{0.33\textwidth}
		\centering
		\includegraphics[width=\textwidth]{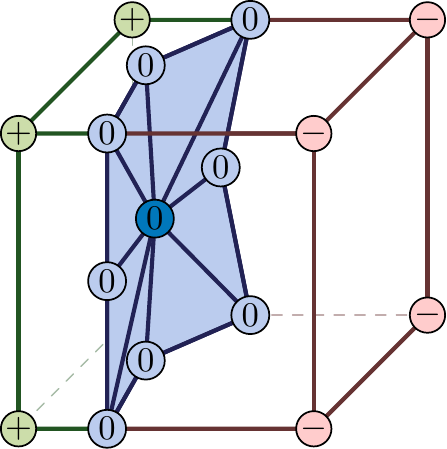}
		\caption{}
		\label{fig:levelset_midpoint_3d}
	\end{subfigure}\hfill%
	\begin{subfigure}[b]{0.33\textwidth}
		\centering
		\includegraphics[width=\textwidth]{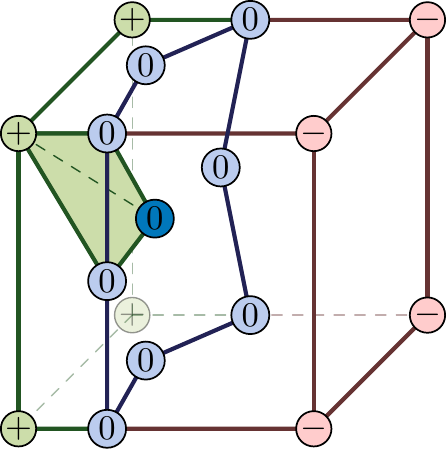}
		\caption{}
		\label{fig:tessellation1_3d}
	\end{subfigure}\hfill%
	\begin{subfigure}[b]{0.33\textwidth}
		\centering
		\includegraphics[width=\textwidth]{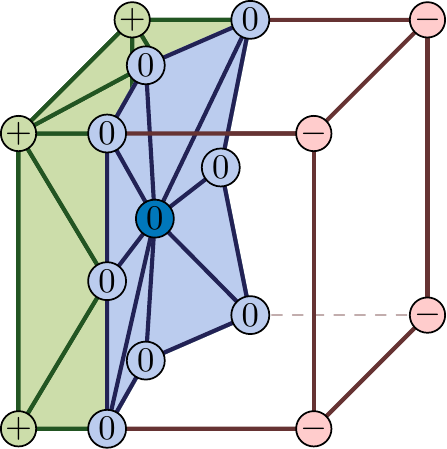}
		\caption{}
		\label{fig:tessellation_3d}
	\end{subfigure}
	\caption{Schematic representation of the midpoint tessellation procedure for a three-dimensional case.}
	\label{fig:midpoint_3d}
\end{figure}

To overcome the deficiencies associated with Delaunay tessellation, we have developed a dedicated tessellation procedure that suits the needs of our immersed analysis approach \cite{divi_error-estimate-based_2020}. We refer to the developed procedure as \emph{midpoint tessellation}, which is illustrated in Figs.~\ref{fig:midpoint_2d} and \ref{fig:midpoint_3d} for the two-dimensional and three-dimensional case, respectively. The \textsc{get\_mosaic\_element} function outlined in Alg.~\ref{alg:midpoint} implements our midpoint tessellation procedure. This function is called by the octree algorithm at the lowest level of subdivision (Alg.~\ref{alg:trim}, L\ref{alg:trim:mosaic}), taking the $2^\ndims$ level set \textit{values} at the octree vertices as input. The function returns a tessellation of the octree cell that (approximately) fits to the immersed boundary.

To explain Alg.~\ref{alg:midpoint}, we first consider the two-dimensional case, which is illustrated in Fig.~\ref{fig:midpoint_2d}. The midpoint tessellation procedure commences with looping over all the edges of the element and recursively calling the \textsc{get\_mosaic\_element} function to truncate the edges that are intersected by the immersed boundary (L\ref{alg:midpoint:startrecursion}-L\ref{alg:midpoint:endrecursion}, Fig.~\ref{fig:trim_edges_2d}). A set of \textit{zero\_points} is then computed by linear interpolation of the level set function across the diagonals between the centroid of the rectangle and its vertices (L\ref{alg:midpoint:zeropoints}, Fig.~\ref{fig:levelset_diagonal_2d}), and the arithmetic average of these points is defined as the \textit{midpoint} (L\ref{alg:midpoint:midpoint}, Fig.~\ref{fig:diagonal_midpoint_2d}). The tessellation is then created by extruding the (truncated) edges toward this \textit{midpoint} (Figs.~\ref{fig:tessellation1_2d}-\ref{fig:tessellation2_2d}). Note that if this procedure is applied to the negated level set values, a tessellation of the complementary part of the rectangular element with a coincident immersed boundary is obtained (Fig.~\ref{fig:tessellation3_2d}).

Since the tessellation algorithm recursively traverses the dimensions of the element, it can directly be extended to the three-dimensional case, as illustrated in Fig.~\ref{fig:midpoint_3d}. In three dimensions, all six faces of the element are tessellated by calling the \textsc{get\_mosaic\_element} function (L\ref{alg:midpoint:startrecursion}-L\ref{alg:midpoint:endrecursion}, Fig.~\ref{fig:diagonal_midpoint_3d}). Based on the diagonals between the centroid of the element and its vertices, zero level set points (L\ref{alg:midpoint:zeropoints}) and a corresponding \textit{midpoint} (L\ref{alg:midpoint:midpoint}) are then computed. The three-dimensional tessellation is finally constructed by extrusion of all (truncated) faces toward the \textit{midpoint} (Fig.~\ref{fig:tessellation1_3d}-\ref{fig:tessellation_3d}). As in the two dimensional case, a conforming interface is obtained when the procedure is applied to the negated level set values.

\begin{remark}[Generalization to non-rectangular elements]
The algorithms presented in this section are presented for the case of hyperrectangles, \emph{i.e.}, a rectangle in two dimensions and a hexagon in three dimensions. The algorithms can be generalized to a broader class of element shapes (\emph{e.g.}, simplices, similar to Ref.~\cite{varduhn_tetrahedral_2016}). For the algorithms to be generalizable, the considered element must be able to define \textit{children} (of the same type). Moreover, the element and its faces must be convex, so that the \textit{midpoint} is guaranteed to be in the interior of the element. Note that this convexity requirement pertains to the shape of the untrimmed element, and not to the trimmed element.
\end{remark}

\subsection{Topology preservation}\label{sec:topology}
\begin{figure}
    \centering
	\begin{subfigure}[b]{0.45\textwidth}
		\centering
	    \includegraphics[width=\textwidth]{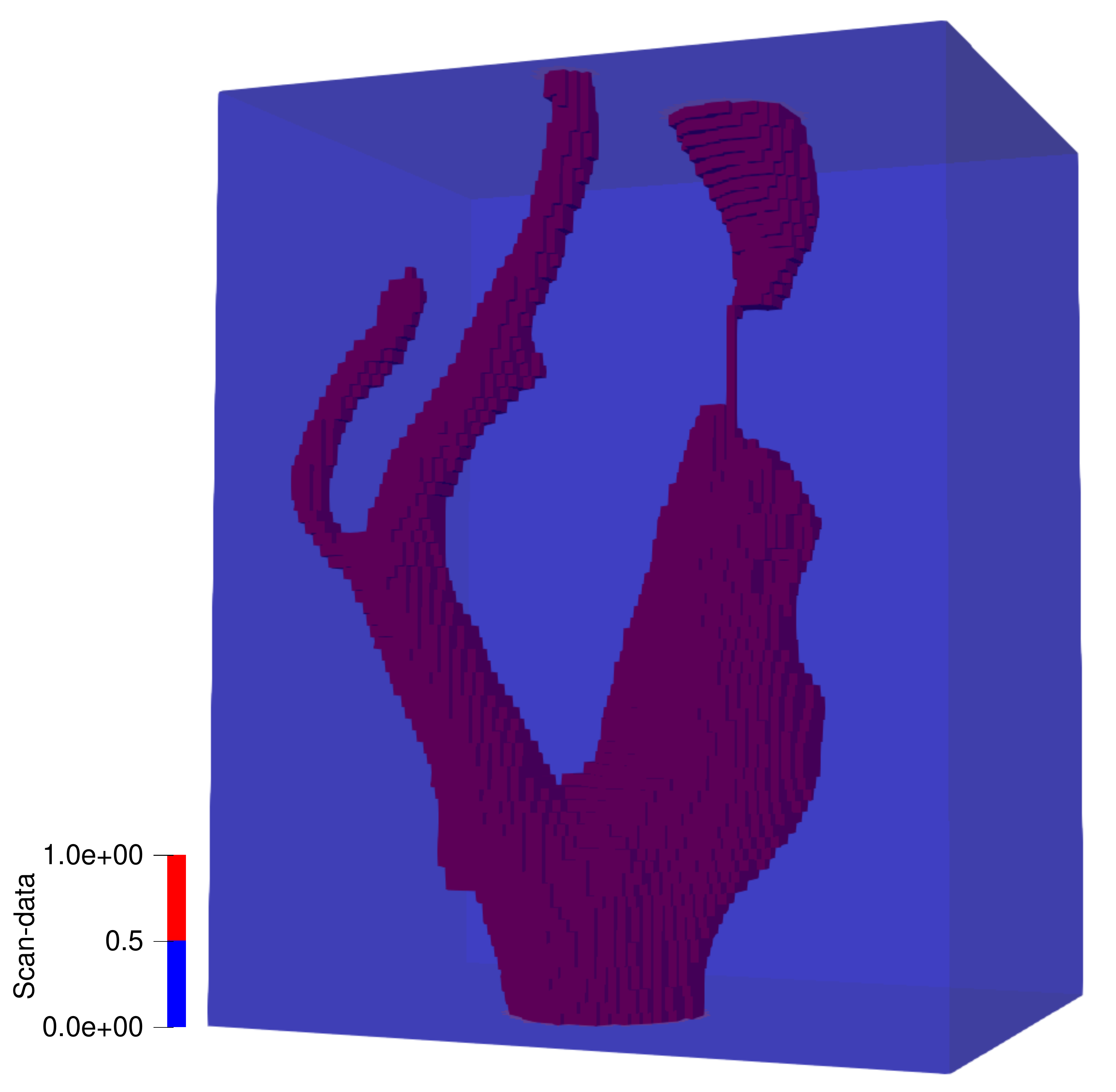}	
		\caption{}
		\label{fig:topology_original}
	\end{subfigure}%
	\hfill
		\begin{subfigure}[b]{0.23\textwidth}
		\centering
		\includegraphics[width=\textwidth]{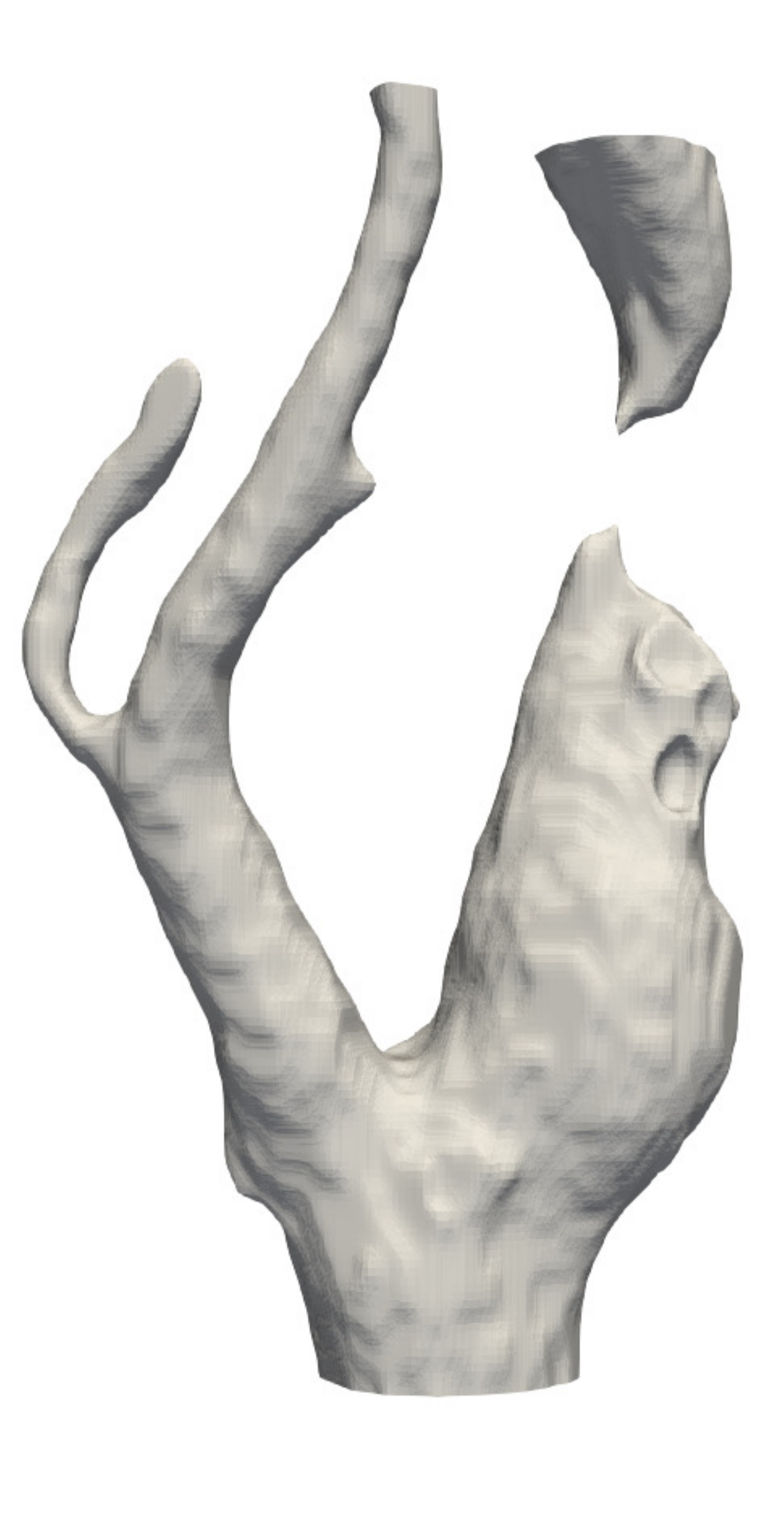}	
		\caption{}
		\label{fig:topology_segmented}
	\end{subfigure}%
	\hfill
		\begin{subfigure}[b]{0.23\textwidth}
		\centering
		\includegraphics[width=\textwidth]{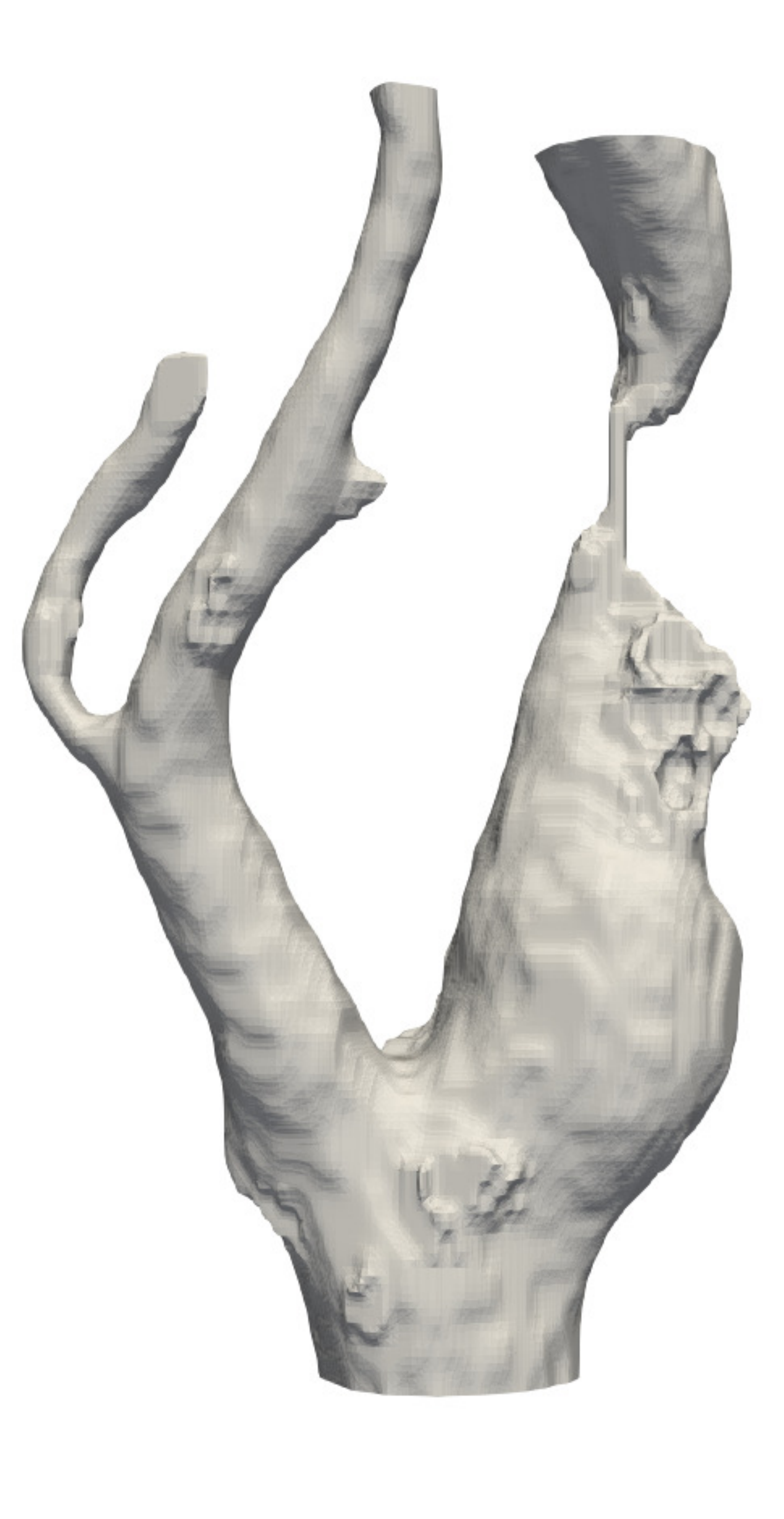}	
		\caption{}
		\label{fig:topology_preserved}
	\end{subfigure}%
    \caption{Example of the spline-based segmentation procedure for a scan-based image of a carotid artery. (\subref{fig:topology_original}) The original voxel data; (\subref{fig:topology_segmented}) the spline segmentation based on the voxel mesh, resulting in a topologically erroneous domain; and (\subref{fig:topology_preserved}) the spline segmentation after application of the topology-preservation algorithm.}
    \label{fig:topology}
\end{figure}

The spline-based segmentation procedure has been demonstrated to yield analysis-suitable domains for a wide range of test cases (see, \emph{e.g.}, Refs.~\cite{verhoosel_image-based_2015,hoang_skeleton-stabilized_2019,de_prenter_condition_2017,de_prenter_preconditioning_2019,divi_error-estimate-based_2020,divi_topology-preserving_2022,divi_residual-based_2022}). An example from Ref.~\cite{divi_topology-preserving_2022} is shown in Fig.~\ref{fig:topology}. This example shows a carotid artery, obtained from a CT-scan. The scan data consists of 80 slices, separated by a distance of $400\,{\rm \mu m}$. Each slice image consists of $85\times 70$ voxels of size $300\times 300 {\rm \mu m}^2$.

When the spline-based segmentation procedure is performed using a B-spline level set function defined on the voxel grid, the result shown in Fig.~\ref{fig:topology_segmented} is obtained. Although the smoothing characteristic of the technique is overall beneficial, in the sense that it leads to smooth boundaries, comparison to the original voxel data in Fig.~\ref{fig:topology_original} shows that in this particular case the topology of the object is altered by the segmentation procedure. This can occur in cases where the features of the object to be described are not significantly larger in size than the voxels (\emph{i.e.}, the Nyquist criterion is not satisfied; see Section~\ref{sec:levelset}). In many cases, the altering of the topology of an object fundamentally changes the problem under consideration, and is therefore generally undesirable.

To avoid the occurrence of topological anomalies due to smoothing, a topology-preservation strategy has been developed in Ref.~\cite{divi_topology-preserving_2022}. The developed strategy follows directly from the smoothing analysis presented in Section~\ref{sec:levelset}, which shows that features with a small relative length scale, $\hat{\ell}=\ell/h$, may be lost upon smoothing, \emph{cf.} equation~\eqref{eq:approxmax}. Hence, topological features may be lost when the mesh size on which the B-spline level set is constructed is relatively large compared to the voxel size. The pivotal idea of the strategy proposed in Ref.~\cite{divi_topology-preserving_2022} is to detect topological anomalies by comparison of the segmented image (Fig.~\ref{fig:topology_segmented}) with the original voxel data (Fig.~\ref{fig:topology_original}) through a moving-window technique. In places where topological anomalies are detected, the mesh on which the smooth level set function is constructed is then refined locally (using THB-splines \cite{giannelli_thb-splines_2012}). This locally increases the relative feature length scale, $\hat{\ell}$, such that the topology is restored (Fig.~\ref{fig:topology_preserved}).

\section{Immersed isogeometric flow analysis}
\label{sec:formulation}

In this section we introduce an immersed discretization of the Stokes flow problem solved on a domain $\domain \subset \mathbb{R}^\ndims$ according to \eqref{eq:implicitdomain}, attained through the scan-based segmentation procedure outlined above. The boundary, $\boundary$, as illustrated in Fig.~\ref{fig:fcmdomains}, is (partly) immersed, in the sense that it does not coincide with element boundaries.

The Stokes flow problem reads
\begin{equation}\label{eq:stokesequations}
\left\{
	\begin{aligned}
		-\nabla \cdot (2\mu \nabla^s \uu ) + \nabla p 	& = \bodystokes  &  &\text{in } \domain, \\
		\nabla \cdot \uu 								& = 0 & &\text{in } \domain, \\
		\uu  		& = \dirichletdatastokes & &\text{on } \dirichletboundary, \\
		2 \mu \left( \nabla^s \uu  \right) \nn - p \nn  			& = \neumanndatastokes & &\text{on } \neumannboundary, 
	\end{aligned}
\right. 
\end{equation}
with velocity $\uu$, pressure $p$, constant viscosity $\mu$, body force $\bodystokes$,  Dirichlet data $\dirichletdatastokes$ and Neumann data $\neumanndatastokes$. The boundary is composed of a Neumann part, $\neumannboundary$, and a Dirichlet part, $\dirichletboundary$, such that  $\overline{\neumannboundary} \cup \overline{\dirichletboundary} = \boundary$ and $\neumannboundary\cap\dirichletboundary=\emptyset$. The vector $\nn$ in the last line denotes the outward-pointing unit normal to the boundary.

When discretizing the Stokes problem \eqref{eq:stokesequations}, the immersed setting poses various challenges: \emph{(i)} Since the Dirichlet (\emph{e.g.}, no-slip) boundary is (partly) immersed, Dirichlet boundary conditions cannot be imposed strongly (\emph{i.e.}, by constraining degrees of freedom) \cite{hansbo_unfitted_2002,embar_imposing_2010}; \emph{(ii)} stability and conditioning issues can occur on unfavorably cut elements \cite{verhoosel_image-based_2015,de_prenter_condition_2017, de_prenter_note_2018, burman_cutfem_2015,de_prenter_stability_2022}; and \emph{(iii)} elements which are known to be inf-sup stable in boundary-fitted finite elements (\emph{e.g.}, Taylor-Hood elements) can lose stability when being cut, resulting in oscillations in the velocity and pressure fields
\cite{hoang_mixed_2017}.

To enable scan-based immersed isogeometric analyses, we have developed a stabilized formulation that addresses these challenges. In this formulation, Dirichlet boundary conditions are imposed weakly through Nitsche's method \cite{hansbo_unfitted_2002,embar_imposing_2010}. Ghost stabilization \cite{burman_ghost_2010} is used to avoid conditioning and stability problems associated with unfavorably cut elements, and skeleton-stabilization is used to avoid inf-sup stability problems. Skeleton-stabilization also allows us to consider equal-order discretizations of the velocity and pressure spaces, simplifying the analysis framework. In Section~\ref{sec:fcmsetting} we first formalize the immersed analysis setting, after which the developed formulation is detailed in Section~\ref{sec:immersogeometric}.

\subsection{Immersed analysis setting}\label{sec:fcmsetting}

The physical domain is immersed in the (cuboid) scan domain, \emph{i.e.}, $\ambientdomain \supset \domain$, on which a locally refined scan mesh $\ambientmesh$ with elements $\element$ is defined. Locally refined meshes can be constructed by sequential bisectioning of selections of elements in the mesh, starting from a Cartesian mesh, which will be discussed in Section~\ref{sec:adaptivesplines}.

\begin{figure}
   \centering
   \begin{subfigure}[b]{0.45\textwidth}
   \centering
   \includegraphics[width=\textwidth]{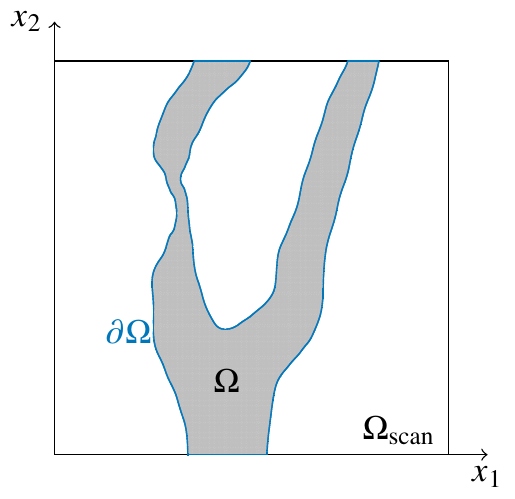}
   \caption{}\label{fig:fcmdomaindef}   
   \end{subfigure}
   \hfill
   \begin{subfigure}[b]{0.45\textwidth}
   \centering
   \includegraphics[width=\textwidth]{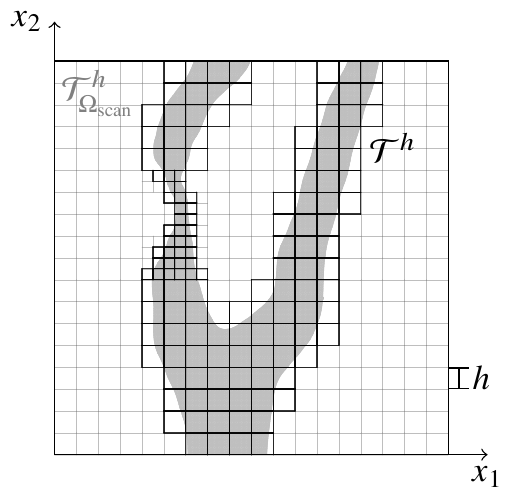}
   \caption{}\label{fig:fcmdomainsa}   
   \end{subfigure}\\[12pt]
   \begin{subfigure}[b]{0.45\textwidth}
   \centering
    \includegraphics[width=\textwidth]{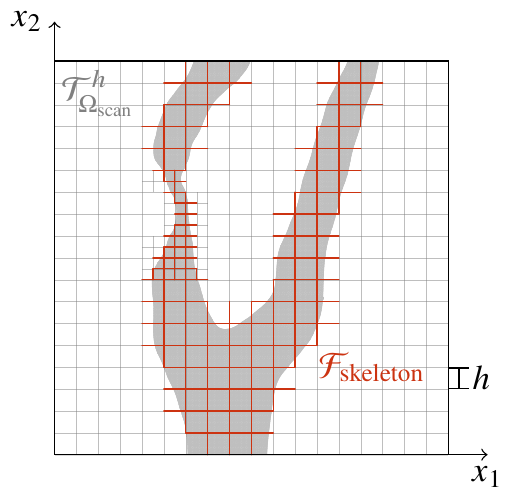}
   \caption{}\label{fig:fcmdomainsb}   
   \end{subfigure}
   \hfill
   \begin{subfigure}[b]{0.45\textwidth}
   \centering
   \includegraphics[width=\textwidth]{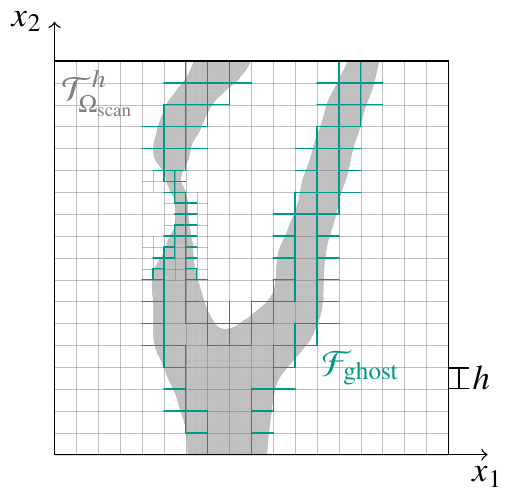}
   \caption{}\label{fig:fcmdomainsc}   
   \end{subfigure}
   
	 \caption{(\subref{fig:fcmdomaindef}) A physical domain $\domain$ (gray), with boundary $\boundary$ (blue), is embedded in the scan domain $\ambientdomain$ (white). (\subref{fig:fcmdomainsa}) The background mesh $\mesh$, which consists of all elements  that intersect the physical domain, is constructed by locally refining the ambient domain mesh $\ambientmesh$. The skeleton mesh, $\mathcal{F}_{\rm skeleton}$, and ghost mesh, $\mathcal{F}_{\rm ghost}$, are shown in panels (\subref{fig:fcmdomainsb}) and (\subref{fig:fcmdomainsc}), respectively.}
     \label{fig:fcmdomains}
\end{figure}

Elements that do not intersect with the physical domain can be omitted, resulting in the locally refined (active) background mesh
\begin{equation}\label{equation:background}
\mesh := \{\element \, | \, \element \in \ambientmesh, \element \cap \domain \neq \emptyset \}.
\end{equation}
The (active) background mesh is illustrated in Fig.~\ref{fig:fcmdomainsa}. By cutting the elements that are intersected by the immersed boundary $\boundary$, a mesh that conforms to the physical domain $\Omega$ is obtained:
\begin{equation}\label{equation:cutmesh}
\cutmesh := \{\element \cap \Omega \, | \, \element \in \mesh \}.
\end{equation}
The tessellation procedure discussed in Section~\ref{sec:tessellation} provides a polygonal approximation of the immersed boundary $\boundary$ through the set of boundary faces 
\begin{equation}\label{equation:boundaryedges}
\boundarymesh := \{ \edge \subset \boundary \, | \, \edge = \partial \element \cap \partial \domain,\,  \element \in \cutmesh \}.
\end{equation}

The considered formulation (see Section \ref{sec:immersogeometric}) incorporates stabilization terms formulated on the edges of the background mesh, which we refer to as the skeleton mesh
\begin{equation}\label{equation:skeletoninterfaces}
\skeleton := \{ F = \partial \element \cap \partial \element' \,|\ \element,\element'\in \mesh, \element\neq \element' \}.
\end{equation}
Note that the faces $F \in \skeleton$ can be partially outside of the domain $\Omega$ and that the boundary of the background mesh is not part of the skeleton mesh. This skeleton mesh is illustrated in Fig.~\ref{fig:fcmdomainsb}.

In addition to the skeleton mesh, we define the ghost mesh as the subset of the skeleton mesh with the faces that belong to  an element intersected by the domain boundary, \emph{i.e.},
\begin{equation}\label{equation:ghostinterfaces}
\ghost := \{ F \cap \partial \element \,|\, F \in \skeleton, \element \in \mathcal{G} \},
\end{equation}
where $\mathcal{G} := \{ K \in \mesh \mid K \cap \partial \Omega \neq \emptyset \}$ is the collection of elements in the background mesh that are crossed by the immersed boundary. The ghost mesh is illustrated in Fig.~\ref{fig:fcmdomainsc}.

\subsection{Stabilized formulation}\label{sec:immersogeometric}
To solve the Stokes problem \eqref{eq:stokesequations} we discretize the velocity and pressure fields using truncated hierarchical B-splines \cite{giannelli_thb-splines_2012,van_brummelen_adaptive_2021}. THB-splines form a basis of degree $\order$ and regularity $\regularity$ constructed over the locally-refined background mesh, $\mesh$, spanning the spline space
\begin{equation}
\splinespace(\mesh) 
= 
\{ \basisfunc\in{}C^{\regularity}(\mesh):\basisfunc|_{\element}\in{}\mathcal{Q}^{\order}(\element),\,\forall{}\element\in \mesh\},
\end{equation}
with $\mathcal{Q}^k(K)$ the set of $d$-variate polynomials on the element $K$ constructed by the tensor-product of univariate polynomials of order $\order$. We consider equal-order discretizations of the velocity and pressure spaces with optimal regularity THB-splines, \emph{i.e.}, $\alpha=k-1$:
\begin{align}
    \uuh &\in \boldsymbol{V}^h = [\mathcal{S}_{k-1}^k]^d \subset [H^1]^d, &   
    \ph \in Q^h &= \mathcal{S}_{k-1}^k \subset L^2.
\end{align}
Note that the superscript $\meshsize$ is used to indicate that these fields are approximations obtained on a mesh with (local) element size $\meshsize$. 

We consider the stabilized Bubnov-Galerkin formulation
\begin{align}
\!\!\!\!\!\!\left\{ \begin{aligned} 
&\text{Find} \: \boldsymbol{u}^h \in \boldsymbol{V}^h \: \text{and} \: p^h \in Q^h \: \text{such that:} \: & &  \\
&a(\boldsymbol{u}^h, \boldsymbol{v}^h) + b(p^h, \boldsymbol{v}^h) + s_{\rm nitsche}(\boldsymbol{u}^h, \boldsymbol{v}^h) + s_{\rm ghost} (\boldsymbol{u}^h, \boldsymbol{v}^h) = f(\boldsymbol{v}^h) & & \forall \boldsymbol{v}^h \in \boldsymbol{V}^h \\
& b(q^h,\boldsymbol{u}^h) - s_{\rm skeleton} (p^h,q^h) = g(q^h) & & \forall q^h \in Q^h
\end{aligned}\right. \label{eq:weak_stokes}
\end{align}
where the bilinear and linear operators are defined as (see Ref.~\cite{hoang_skeleton-stabilized_2019} for details)
\begin{subequations}
\begin{align}
a(\boldsymbol{u}^h,\boldsymbol{v}^h) & := 2 \mu (\nabla^s \uuh, \nabla^s \vvh)  \nonumber \\
&\phantom{:=} -  2 \mu \left[ \langle (\nabla^s \uuh ) \nn, \vvh \rangle_{\dirichletboundary} + \langle (\nabla^s \vvh ) \nn, \uuh \rangle_{\dirichletboundary} \right]\\
b(p^h,\boldsymbol{v}^h) &:= - (\ph, \nabla \cdot \vvh) + \langle p^h, \vvh \cdot \nn \rangle_{\dirichletboundary}\\
f (\boldsymbol{v}^h) &:= (\bodystokes, \vvh) + \langle \neumanndatastokes, \vvh \rangle_{\neumannboundary} - 2 \mu \langle  (\nabla^s \vvh) \nn, \dirichletdatastokes \rangle_{\dirichletboundary} \\
g (q^h) &:= \langle \qh, \dirichletdatastokes \cdot \nn \rangle_{\dirichletboundary} \\
s_{\rm nitsche} (\boldsymbol{u}^h,\boldsymbol{v}^h)  &:=  \langle \beta \mu  \meshsize^{-1} ( \uuh - \dirichletdatastokes ) ,\vvh \rangle_{\dirichletboundary}  \label{eq:nitscheoperator}\\
s_{\rm ghost} (\boldsymbol{u}^h,\boldsymbol{v}^h) &:= \sum \limits_{F \in \mathcal{F}_{\rm ghost}} \int_{F} \gamma_g \mu h_F^{2k-1} \llbracket \partial_n^k \boldsymbol{u}^h \rrbracket \cdot \llbracket \partial_n^k \boldsymbol{v}^h \rrbracket \: {\rm d}S\label{eq:ghostoperator}\\
s_{\rm skeleton} (p^h,q^h)  &:= \sum_{F \in \mathcal{F}_{\rm skeleton}} \int_{F} \gamma_{s} \mu^{-1} h_F^{2k+1} \llbracket \partial_n^k p^h \rrbracket \llbracket \partial_n^k q^h \rrbracket \: {\rm d}S \label{eq:skeletonoperator}
\end{align}%
\label{eq:operators}%
\end{subequations}
where $(\cdot,\cdot)$ denotes the inner product in $L^{2}(\domain)$,  $\langle \cdot, \cdot \rangle_{\boundary}$ denotes the inner product in $L^{2}(\boundary)$, and $\llbracket \cdot \rrbracket$ denotes the interface jump operator. The parameters $\beta$, $\gamma_g$, and $\gamma_s$ denote the penalty constants for the Nitsche term, the ghost-stabilization term, and the skeleton-stabilization term, respectively. 

To ensure stability and optimal approximation, the Nitsche stabilization term \eqref{eq:nitscheoperator} scales with the inverse of the (background) mesh size parameter, $\meshsize$ \cite{burman_cutfem_2015}. The Nitsche stability parameter $\beta$ should be selected appropriately, being large enough to ensure stability, while not being too large to cause locking-type effects (see, \emph{e.g.}, Refs.~\cite{de_prenter_note_2018,badia_mixed_2018,de_prenter_stability_2022}). The ghost-penalty operator in \eqref{eq:ghostoperator} controls the $k^{\rm th}$-order normal derivative jumps over the interfaces of the elements which are intersected by the domain boundary $\boundary$. Since in this contribution splines of degree $\order$ with $C^{\order-1}$-continuity are considered, only the jump in the $\order^{\rm th}$ normal derivative is non-vanishing at the ghost mesh. The ghost-stabilization term scales with the size of the faces as $\meshsize_F^{2k-1}$. Appropriate selection of the parameter $\gamma_{g}$ corresponding with the Nitsche parameter, $\beta$, assures the stability of the formulation independent of the cut-cell configurations. To avoid loss of accuracy, the ghost-penalty parameter, $\gamma_{g}$, should also not be too large \cite{badia_linking_2022}.

The skeleton-stabilization operator \eqref{eq:skeletonoperator}, proposed in Ref.~\cite{hoang_skeleton-stabilized_2019}, penalizes jumps in higher-order pressure gradients. This ensures inf-sup stability of the equal-order velocity-pressure discretization, and resolves spurious pressure oscillations caused by cut elements. This spline-based skeleton-stabilization technique can be regarded as the higher-order continuous version of the interior penalty method proposed in Ref.~\cite{burman_edge_2006}. To ensure stability and optimality, the operator \eqref{eq:skeletonoperator} scales with $\meshsize_F^{2k+1}$. The parameter $\gamma_{s}$ should be selected such that oscillations are suppressed, while the influence of the additional term on the accuracy of the solution remains limited. It is noted that since the inf-sup stability problem is not restricted to the immersed boundary, the skeleton stabilization pertains to all interfaces of the background mesh.

In our scan-based analysis workflow it is, from a computational effort point of view, generally impractical to evaluate the (integral) operators \eqref{eq:operators} exactly. The error of the Galerkin solution with inexact integration, $u^h_{\mathcal{Q}}=(\boldsymbol{u}^h_{\mathcal{Q}},p^h_{\mathcal{Q}})$, is then composed of two parts, \emph{viz.}: \emph{(i)} the discretization error, defined as the difference between the analytical solution, $u =(\boldsymbol{u},p)$, and the approximate Galerkin solution in the absence of integration errors, $u^h =(\boldsymbol{u}^h,p^h)$; and \emph{(ii)} the inconsistency error related to the integration procedure, which is defined as the difference between the approximate solution in the absence of integration errors, $u^h$, and the Galerkin solution with integration errors, ${u}^{h}_{\mathcal{Q}}$. In practice, one needs to control both these error contributions in order to ensure the accuracy of a simulation result. From the perspective of computational effort, it is in general not optimal to make either one of the contributions significantly smaller than the other.

The decoupling of the geometry description from the analysis mesh provides the immersed (isogeometric) analysis framework with the flexibility to locally adapt the resolution of the solution approximation without the need to reparametrize the domain. To leverage this flexibility in the scan-based analysis setting, it is essential to automate the adaptivity procedure, as manual selection of adaptive cut-element quadrature rules and mesh refinement regions is generally impractical on account of the complex volumetric domains that are considered.

In our work we have developed error-estimation-based criteria that enable adaptive scan-based analyses. In Section~\ref{sec:adaptivequadrature} we first discuss an adaptive octree quadrature procedure used to reduce the computational cost associated with cut element integration. In Section~\ref{sec:adaptivesplines} we then discuss a residual-based error estimator to refine the THB-spline approximation of the field variables only in places where this results in substantial accuracy improvements.

\section{Adaptive integration of cut elements}\label{sec:adaptivequadrature}

From the perspective of computational effort, a prominent challenge in immersed finite element methods is the integration of the cut elements. While quadrature points can be constructed directly on all octree sub-cells (Section~\ref{sec:octree}), this generally results in very expensive integration schemes, especially for three-dimensional problems \cite{divi_error-estimate-based_2020}. A myriad of techniques have been developed to make cut-element integration more efficient, an overview of which is presented in, \emph{e.g.}, Refs.~\cite{abedian_performance_2013,divi_error-estimate-based_2020}. In the selection of an appropriate cut element integration scheme one balances robustness (with respect to cut element configurations), accuracy, and expense.

In the context of scan-based analyses, we have found it most suitable to leverage the robustness of the octree procedure as much as possible. To improve the computational efficiency of the resulting quadrature rules, we have developed a procedure that adapts the number of integration points on each integration sub-cell, similar to the approach used in Ref.~\cite{abedian_finite_2013}, as lower-order integration on very small sub-cells does not significantly reduce the accuracy.

\subsection{Integration error estimate}
The pivotal idea of our adaptive octree quadrature procedure is to optimally distribute integration points over the sub-cells using an error estimator based on Strang's first lemma \cite{strang_analysis_2008} (see also \cite[Lemma 2.27]{ern_theory_2004}). In the immersed analysis setting, this lemma provides an upper bound for the error $u-u^h_{\mathcal{Q}}$. Following the derivation in Ref.~\cite{divi_error-estimate-based_2020} (to which we refer the reader for details), this error bound can be expressed in abstract form as
\begin{equation}
\begin{aligned}
 \left\| u - {u}^h_{\mathcal{Q}} \right\|_{W(h)} 
 \leq  
  \vphantom{\sum \limits_{K \in \mathcal{T}^h_\Omega}} \left( 1 + \frac{\left\| a^h \right\|_{W(h),V^h}}{\alpha^h}  \right) \left\| u - \mathcal{I}^h u \right\|_{W(h)} +  
   \frac{1}{\alpha^h} 
\sum \limits_{K \in \mathcal{T}^h_\Omega} \left(e_K^a + e_K^f \right),
\end{aligned}
\label{eq:strang1estimatefinal}
\end{equation}
where $\alpha^h$ denotes the inf-sup constant associated with the (aggregate) bilinear form $a^h:W^h \times V^h \rightarrow \mathbb{R}$, with trial and test velocity-pressure spaces $W^h$ and $V^h$, with $W(h)=\operatorname{span}{}\{u\}\oplus{}W^h$ the linear space containing the weak solution, $u$. The element-integration-error indicators associated with the (aggregate) bilinear form $a^h$ and (aggregate) linear form $f^h$ are respectively elaborated as
\begin{subequations}
\begin{align}
e_K^a &= 
\sup_{v^h_K \in {V}^h_K}{\frac{ \left|\int_{K} A^h_\Omega(\mathcal{I}^h u,v^h_K)(\boldsymbol{x}_K)\,{\rm d}V -  \sum \limits_{l=1}^{l_K} \omega^l_K A^{h}_\Omega(\mathcal{I}^h u,v^h_K )( \boldsymbol{x}_K^l)\right| }{ \left\| v^h_K \right\|_{V^h_K}}}, \label{eq:elem_interrora}\\
e_K^f &= \sup_{v^h_K \in {V}^h_K}{\frac{  \left|\int_{K} F^h_\Omega( v^h_K)(\boldsymbol{x})\,{\rm d}V -  \sum \limits_{l=1}^{l_K} \omega^l_K F^{h}_\Omega(v^h_K)( \boldsymbol{x}_K^l) \right| }{ \left\| v^h_K \right\|_{V^h_K}}}, 
\label{eq:elem_interrorb}
\end{align}
\label{eq:elem_interror}%
\end{subequations}
where $A^{h}_\Omega$ and $F^{h}_\Omega$ are the integrands corresponding to the volumetric terms in the bilinear form and linear form in the Galerkin problem, respectively,  and where for each element $K$, the set $\{(\boldsymbol{x}_K^l,\omega_K^l)\}_{l=1}^{l_K}$ represents a quadrature rule. The norm $\| \cdot \|_{V_K^h}$ corresponds to the restriction of the space $V^h$ to the element $K$. We note that it has been assumed that integration errors associated with the boundary terms in the Galerkin problem are negligible in comparison to the errors in the volumetric terms, which is in line with the goal of optimizing the volumetric quadrature rules of cut elements.

It is desirable to apply a single integration scheme for all terms in the bilinear and linear forms and, hence, to treat the element-integration errors \eqref{eq:elem_interror} in the same way. To do this, we note that the integrals between the absolute bars in the numerators of \eqref{eq:elem_interror} constitute linear functionals on $V^h_{K}$. By the Riesz-representation theorem \cite{ern_theory_2004}, there exist functions $T^{a}, T^{f} \in V^h_{K}$ such that
\begin{align}
\int \limits_{K} T^{a} v^h_K \,{\rm d}V &= \int \limits_{K} A^h_\Omega\left( \mathcal{I}^h u, v^h_K \right) \,{\rm d}V,  & 
\int \limits_{K} T^{f} v^h_K \,{\rm d}V &= \int \limits_{K} F^h_\Omega\left( v^h_K \right) \,{\rm d}V, \label{eq:riesz_representation}%
\end{align}
for all $v^h_K \in V^h_K$. Assuming that the difference in applying the integral quadrature to the left- and right-hand-side members of \eqref{eq:riesz_representation} is negligible, it then holds that
\begin{subequations}
\begin{align}
e_K^{a} & \leq \| T^a \|_{L^2(K)} \sup_{T^a,v^h_K \in V^h_K}{ \frac{\left|\int_{K} T^{a} v^h_K \,{\rm d}V -  \sum \limits_{l=1}^{l_K} \omega^l_K \left(T^{a} v^h_K \right) (\boldsymbol{x}_K^l)\right| }{ \| T^a \|_{L^2(K)} \| v^h_K \|_{V^h_K}}}, \label{eq:interror_riesz_a}\\
e_K^{f} & \leq \| T^f \|_{L^2(K)} \sup_{T^f,v^h_K \in V^h_K}{ \frac{\left|\int_{K} T^{f} v^h_K \,{\rm d}V -  \sum \limits_{l=1}^{l_K} \omega^l_K \left(T^{f} v^h_K \right) (\boldsymbol{x}_K^l)\right| }{ \| T^f \|_{L^2(K)} \| v^h_K \|_{V^h_K}}}. \label{eq:interror_riesz_b}
\end{align}%
\label{eq:interror_riesz}%
\end{subequations}
With both $T^a$ (resp.~$T^f$) and $v^h_K$ in the polynomial space $V^h_K$, the product $T^a v^h_K$ (resp.~$T^f v^h_K$) resides in the double-degree (normed) polynomial space $\mathcal{Q}_K^{2k}$. It then follows that 
\begin{align}
e_K^a &\lesssim \| T^a \|_{L^2(K)} e_K^p,  & e_K^f &\lesssim \| T^f \|_{L^2(K)} e_K^p,
\end{align}
with the uniformly applicable polynomial integration error defined as
\begin{align}
e_K^{p} &= \sup_{p_K \in \mathcal{Q}^{2\order}_K}{\frac{  \left|\int_{K} p_K(\boldsymbol{x}_K)\,{\rm d}V -  \sum \limits_{l=1}^{l_K} \omega^l_K {p_K}(\boldsymbol{x}_K^l)\right| }{ \| p_K \|_{P_K} } } \nonumber \\  &= \left|\int_{K} p_{K,{\rm max}}(\boldsymbol{x}_K)\,{\rm d}V -  \sum \limits_{l=1}^{l_K} \omega^l_K {p_{K,{\rm max}}}(\boldsymbol{x}_K^l)\right| . \label{eq:polyerror}
\end{align}
The supremizer can be evaluated in terms of a polynomial basis $\boldsymbol{\Phi}$ for the space $\mathcal{Q}_K^{2\order}$ as (see Ref.~\cite{divi_error-estimate-based_2020} for a detailed derivation)
\begin{align}
  p_{K,{\rm max}} =  \displaystyle \frac{\mathbf{\Phi}^T  \mathbf{G}^{-1} ( \boldsymbol{\xi} - \bar{\boldsymbol{\xi}}  )}{  \| \boldsymbol{\xi} - \bar{\boldsymbol{\xi}} \|_{\mathbf{G}^{-1}} }, \label{eq:evaluatedpolynomial}
\end{align}
where $\boldsymbol{\xi}=\int_{K} \boldsymbol{\Phi}\,{\rm d}V$, $\bar{\boldsymbol{\xi}}=\sum_{l=1}^{l_K} \omega^l_K {\boldsymbol{\Phi}}(\boldsymbol{x}_K^l)$, and $\mathbf{G}$ is the (positive definite) Gramian matrix associated with the inner product with which the polynomial space is equipped.

\subsection{Quadrature optimization algorithm}
The computable error definition \eqref{eq:polyerror} and the corresponding computable "worst possible" function \eqref{eq:evaluatedpolynomial} form the basis of our adaptive integration procedure, which we summarize in Alg.~\ref{alg:quadrature} (a detailed version is presented in Ref.~\cite{divi_error-estimate-based_2020}). The developed optimization procedure is intended as a per-element preprocessing operation, which results in optimized \textit{quadrature} rules for all cut elements in a mesh. Besides the \textit{partition}, $\mathcal{P}_K$, of element, $K$, the procedure takes the \textit{order} of the monomial basis, $\mathcal{Q}_K$, and stopping \textit{criterion} (\emph{e.g.}, a prescribed number of integration points) as input.

\begin{algorithm}
\caption{Function to optimize the distribution of cut element quadrature points}\label{alg:quadrature}
    \hspace*{\algorithmicindent} \textbf{Input:} element \textit{partition}, basis function order \textit{order}, stopping \textit{criterion}\\
    \hspace*{\algorithmicindent} \textbf{Output:} optimized \textit{quadrature} rule
\begin{algorithmic}[1]
\Function{optimize\_quadrature}{\textit{partition}, \textit{order}, \textit{criterion}}
\State \textit{basis} = \Call{get\_monomial\_basis}{\textit{order}}\label{alg:quadrature:initbasis}
\State \textit{xi\_exact} = \Call{exact\_integration}{\textit{basis}, \textit{order}}\label{alg:quadrature:xiexact}
\State \textit{gramian} = \Call{get\_gramian\_matrix}{\textit{basis}, \textit{order}}\label{alg:quadrature:gramian}
\State \textit{quadrature} = \Call{initialize\_quadrature}{\textit{partition}}\label{alg:quadrature:initquadrature}
\While{\textbf{not} \textit{criterion}}\label{alg:quadrature:whilecriterion}
\Comment{Adapt quadrature until the stopping criterion is met}
    \State \textit{xi\_quadrature} = \Call{quadrature\_integration}{\textit{basis}, \textit{quadrature}}\label{alg:quadrature:xiquadrature}
    \State \textit{worst\_function} = \Call{get\_worst\_function}{\textit{xi\_exact}, \textit{xi\_quadrature}, \textit{gramian}, \textit{basis}}\label{alg:quadrature:worstfunction}
    \State \textit{indicators} = \Call{initialize\_indicators}{\textit{partition}}
    \For{\textit{subcell}, \textit{indicator} \textbf{in} \Call{zip}{\textit{partition}, \textit{indicators}}}
    \Comment{Iterate over the sub-cells}\label{alg:quadrature:subcellloop}
        \State \textit{error} = \Call{get\_subcell\_error}{\textit{worst\_function}, \textit{quadrature}}\label{alg:quadrature:error}
        \State \textit{cost} = \Call{get\_subcell\_cost}{\textit{quadrature}}\label{alg:quadrature:cost}
        \State \textit{indicator} = \Call{get\_subcell\_indicator}{\textit{subcell\_error}, \textit{subcell\_cost}}\label{alg:quadrature:indicators}
    \EndFor
    \State \textit{marking} = \Call{mark\_subcells}{\textit{indicators}}
    \Comment{Mark based on marking strategy}\label{alg:quadrature:marking}
    \State \textit{quadrature} = \Call{update\_quadrature}{\textit{quadrature}, \textit{marking}}\label{alg:quadrature:update}
\EndWhile
\State \textbf{return} \textit{quadrature}\label{alg:quadrature:return}
\EndFunction
\end{algorithmic}
\end{algorithm}

The procedure commences with the determination of the polynomial \textit{basis}, $\mathbf{\Phi}$ (L\ref{alg:quadrature:initbasis}), the evaluation of the basis function integrals, $\boldsymbol{\xi}$ (L\ref{alg:quadrature:xiexact}), the computation of the \textit{gramian} matrix, $\mathbf{G}$ (L\ref{alg:quadrature:gramian}), and the initialization of the partition \textit{quadrature} rule (L\ref{alg:quadrature:initquadrature}). This initial quadrature rule corresponds to the case where the lowest order (one point) integration rule is used on each sub-cell in the partition. It is noted that the integral computations with full Gaussian quadrature for the \textit{basis} and \textit{gramian} are relatively expensive, but that the computational efficiency gains from the optimized integration scheme outweigh these costs when used multiple times.

The error-estimation-based quadrature optimization is then performed in an incremental fashion (L\ref{alg:quadrature:whilecriterion}), until the stopping \textit{criterion} is met. Given the considered \textit{quadrature} rule, the approximate basis integrals (L\ref{alg:quadrature:xiquadrature}) and worst possible function to integrate \eqref{eq:evaluatedpolynomial} (L\ref{alg:quadrature:worstfunction}) are determined. Subsequently, for each sub-cell, $\wp$, in the \textit{partition} (L\ref{alg:quadrature:subcellloop}), on L\ref{alg:quadrature:error} the sub-cell error indicators 
\begin{align}
e_\wp^{p}  &= \left|\int_{\wp} p_{K,{\rm max}}(\boldsymbol{x}_K)\,{\rm d}V -  \sum \limits_{l \in \mathcal{I}_{\wp}} \omega^l_K {p_{K,{\rm max}}}(\boldsymbol{x}_K^l)\right| 
\end{align}
are computed. In this expression, $\mathcal{I}_{\wp}$, is the set of indices corresponding to integration points on the sub-cell $\wp$. Note that the sum of the sub-cell errors, $e_\wp^{p}$, provides an upper bound for the element integration error \eqref{eq:polyerror}. Sub-cell \textit{indicators} are then computed (L\ref{alg:quadrature:indicators}) by weighing the sub-cell errors with the \textit{costs} associated with increasing the quadrature order on a particular sub-cell, as evaluated on L\ref{alg:quadrature:cost} (see Ref.~\cite{divi_error-estimate-based_2020} for details).

Once the \textit{indicators} have been computed for all sub-cells in the \textit{partition}, the sub-cells with the largest indicators are marked for increasing the number of integration points (L\ref{alg:quadrature:marking}). We propose two marking strategies, \emph{viz.} a sub-cell marking strategy in which only the sub-cell with the highest indicator is marked, and an octree-level marking strategy in which all sub-cells in the octree level with the highest error are marked. After marking, the quadrature order on the marked sub-cells is increased (L\ref{alg:quadrature:update}).

\subsection{Optimized quadrature results}
\begin{figure}
	\centering
	\begin{subfigure}[b]{0.42\textwidth}
		\centering
		\includegraphics[width=0.8\textwidth]{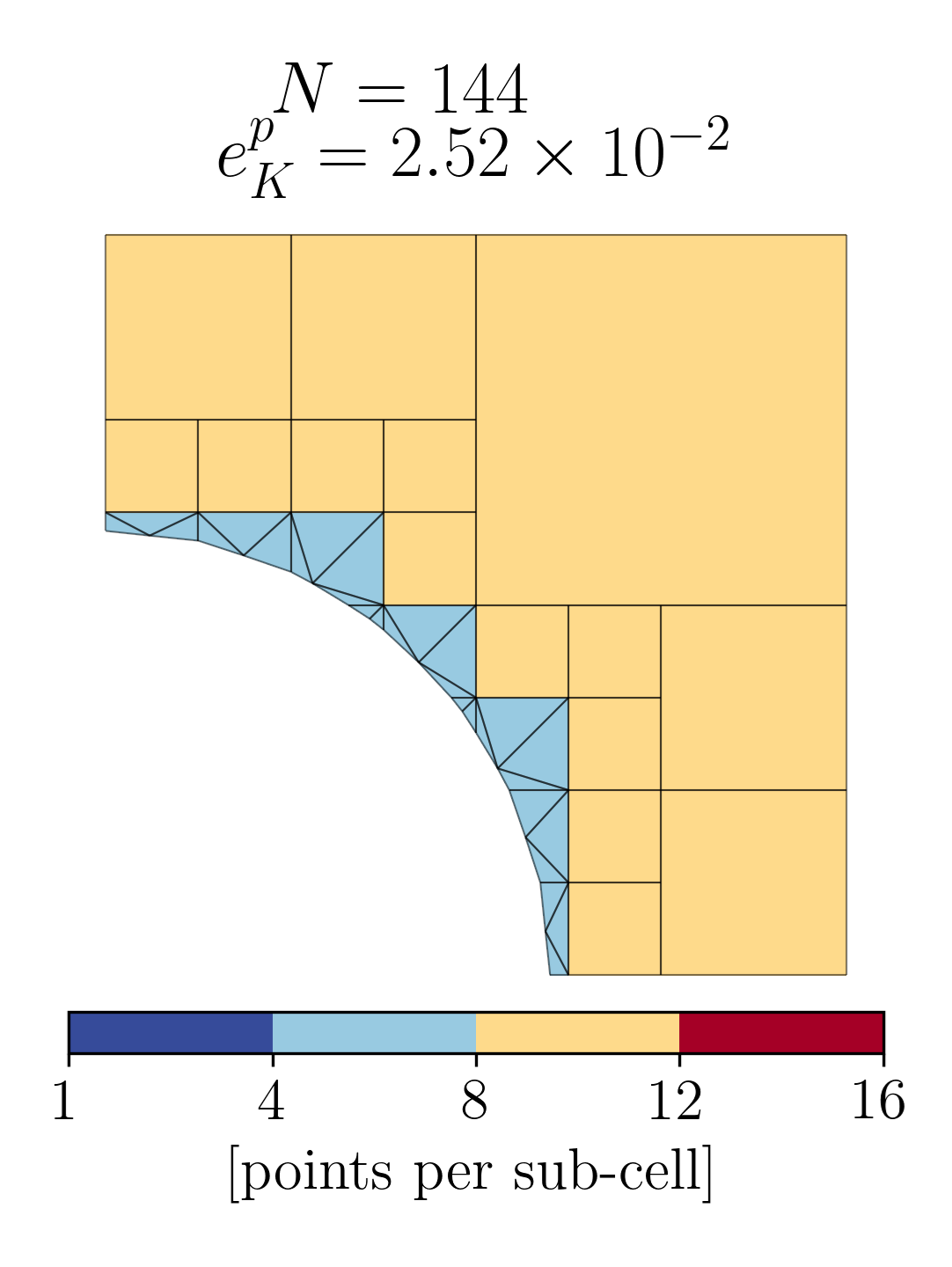}
		\caption{Equal-order Gauss}
		\label{fig:equalorder_gauss}
	\end{subfigure}\hfill%
	\begin{subfigure}[b]{0.42\textwidth}
		\centering
		\includegraphics[width=0.8\textwidth]{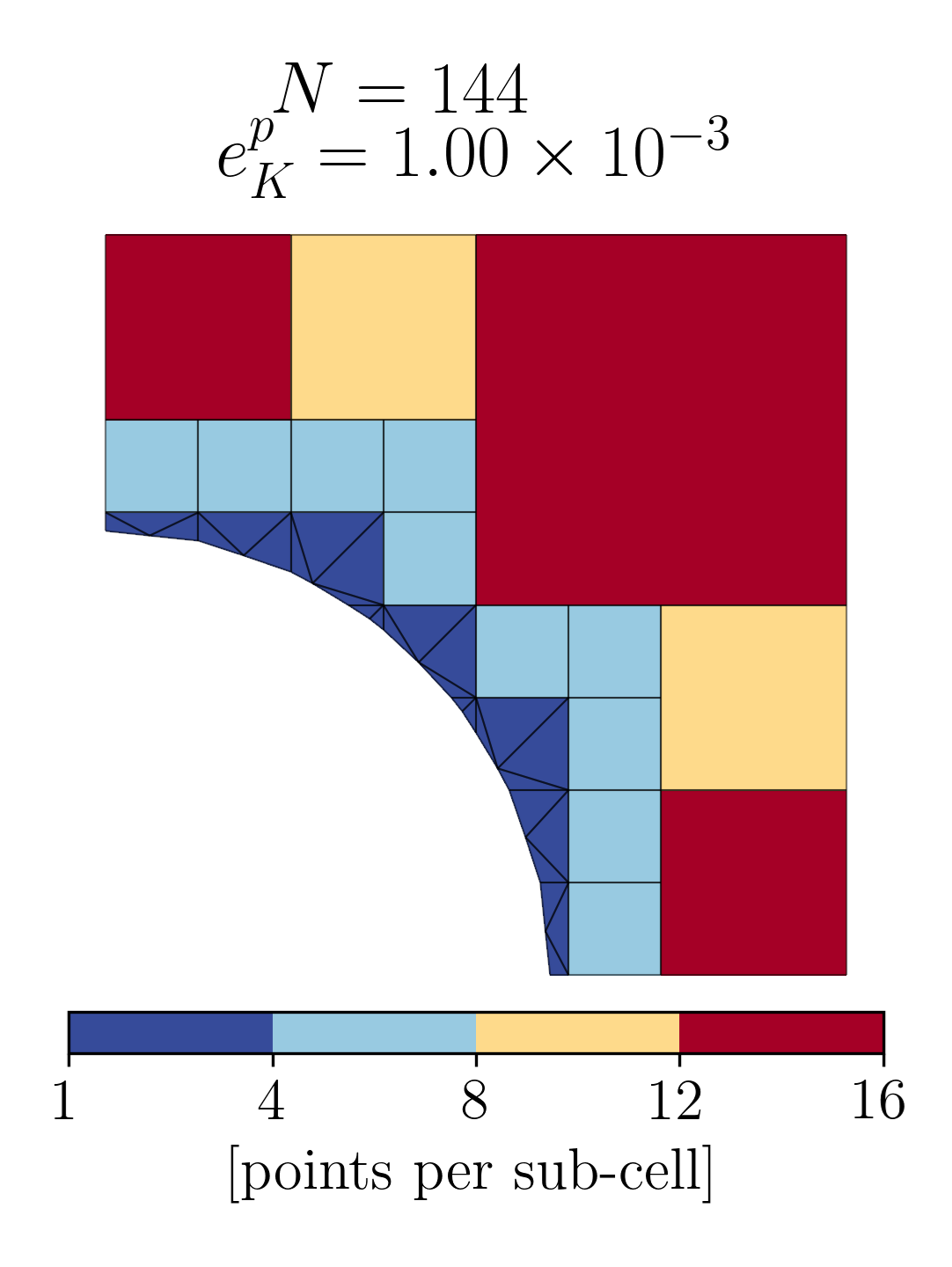}
		\caption{Optimal integration}
		\label{fig:optimal_gauss}
	\end{subfigure}
	\caption{Distribution of integration points over a cut element with $144$ points, comparing the case of (\subref{fig:equalorder_gauss}) an equal-order Gauss scheme, and (\subref{fig:optimal_gauss}) optimally distributed Gauss points using the sub-cell marking strategy. Note that the error is reduced by a factor of 25 by using the sub-cell marking strategy.}
	\label{fig:optimaldist2d}
\end{figure}

A detailed study of the error-estimation-based quadrature optimization scheme is presented in Ref.~\cite{divi_error-estimate-based_2020}. We here reproduce a typical result for a unit square with a circular exclusion, as illustrated in Fig.~\ref{fig:optimaldist2d}. To assess the performance of the developed adaptive integration technique, we study its behavior in terms of integration accuracy versus the number of integration points. 

\begin{figure}
	\centering
	\begin{subfigure}[b]{0.5\textwidth}
		\includegraphics[width=\textwidth]{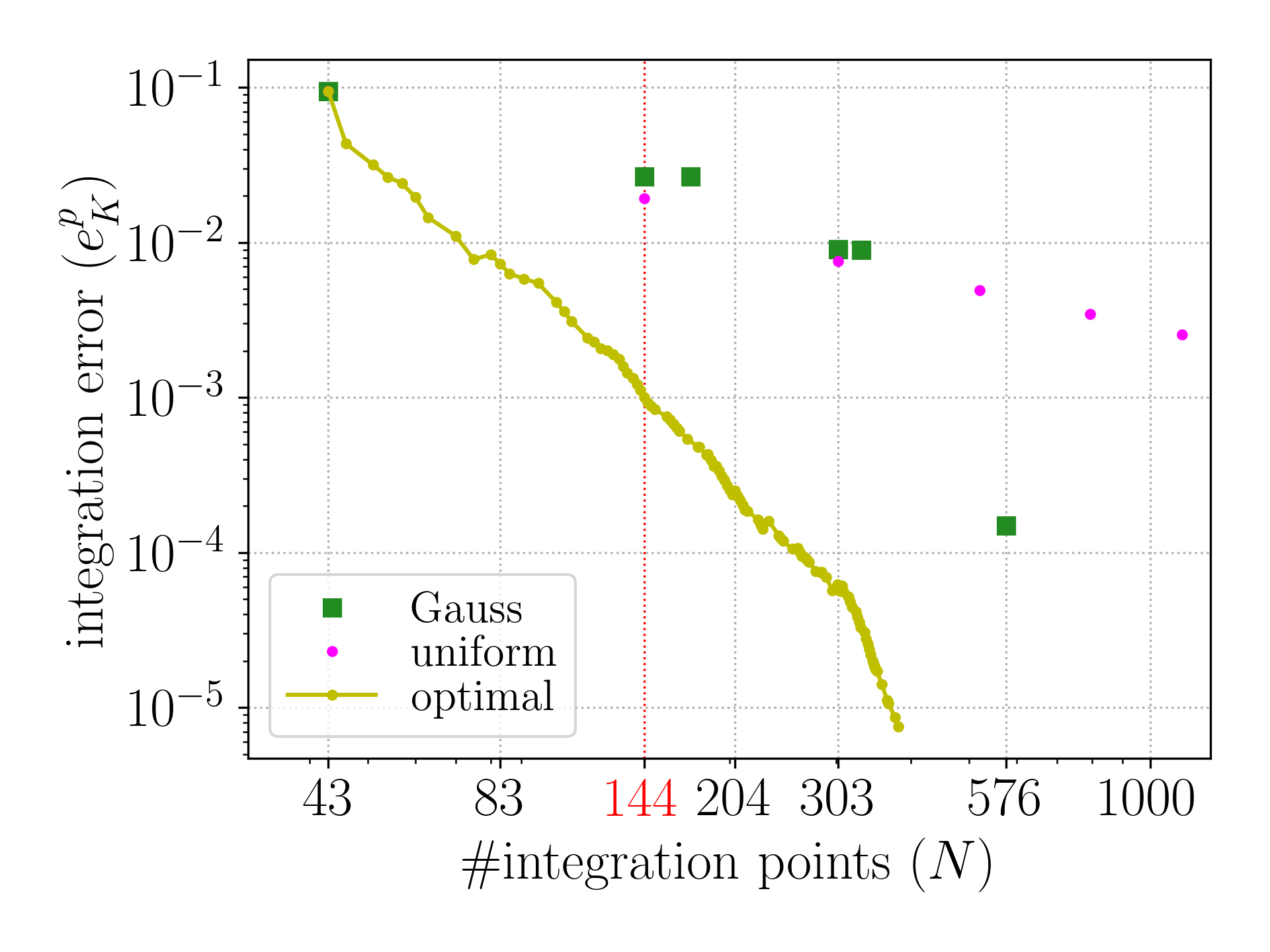}
		\caption{}
		\label{fig:errorvspointsa}
	\end{subfigure}\hfill%
	\begin{subfigure}[b]{0.5\textwidth}
		\includegraphics[width=\textwidth]{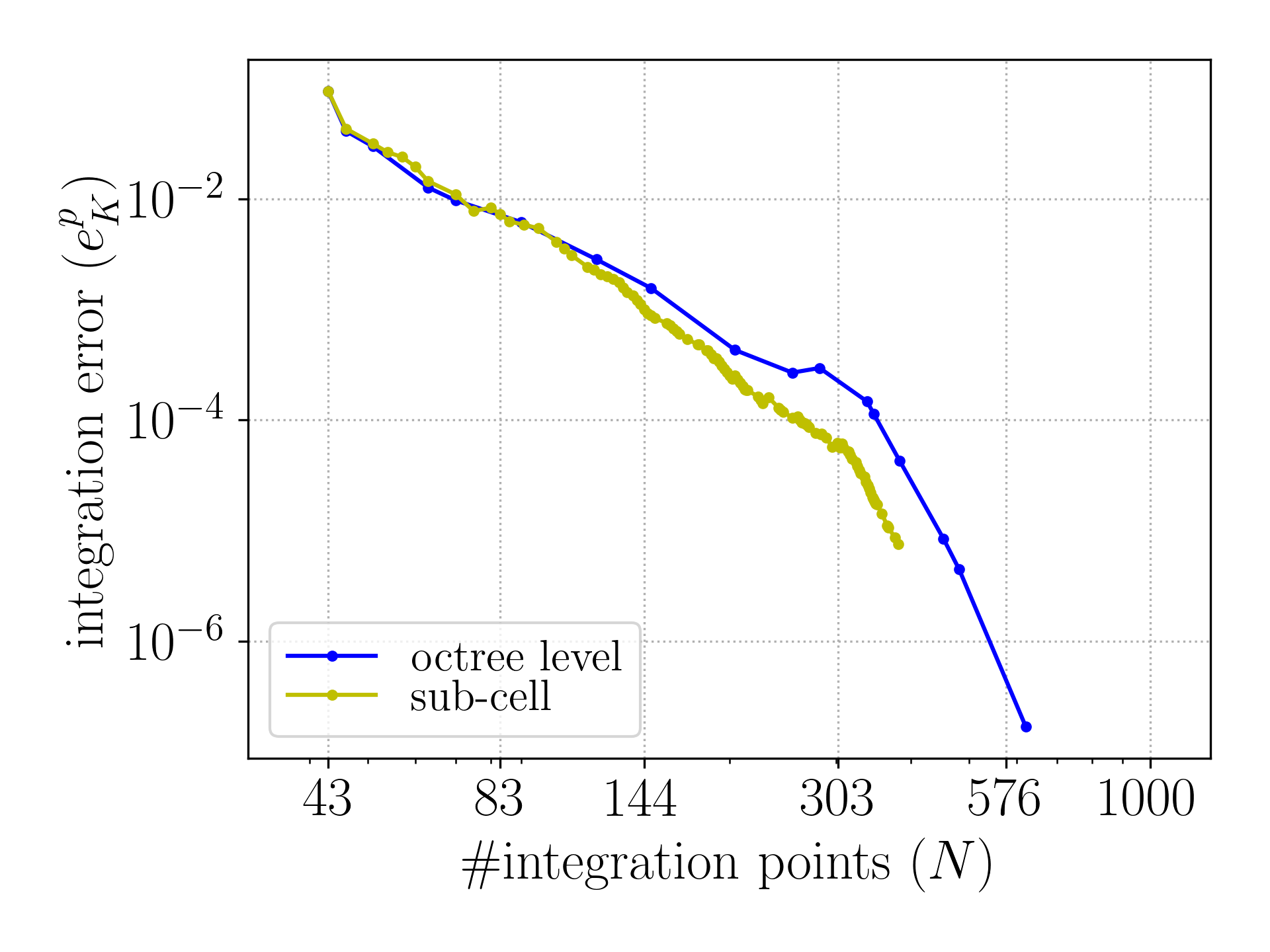}
		\caption{}
		\label{fig:errorvspointsb}
	\end{subfigure}
 	\caption{Integration error \emph{vs.} the number of integration points. (\subref{fig:errorvspointsa}) Comparison of the optimized quadrature results with (non-optimized) equal-order integration. (\subref{fig:errorvspointsb}) Comparison of the sub-cell and octree-level marking strategies.}
 	\label{fig:errorvspoints}
\end{figure}

Fig.~\ref{fig:errorvspointsa} displays the integration error as evolving during the optimization procedure when using the sub-cell marking strategy. The non-optimized case in which the same integration scheme is considered on each sub-cell is displayed for reference. As can be seen, the error associated to the same number of integration points is substantially lowered using the adaptive integration procedure. For example, for the case where $144$ integration points are considered, the error corresponding to the non-optimized second-order Gauss scheme is equal to $2.52 \times 10^{-2}$, while the error corresponding to the optimized quadrature is equal to $1.00 \times 10^{-3}$, \emph{i.e.}, a factor 25 reduction in error. Fig.~\ref{fig:optimaldist2d} displays the distribution of the integration points over the sub-cells for the equal-order Gauss scheme and the optimized case, which clearly demonstrates that the significant reduction in error is achieved by assigning more integration points to the larger sub-cells before introducing additional points in the smaller sub-cells. From Fig.~\ref{fig:errorvspointsa} it is also observed that when the optimization algorithm is terminated at a fixed error of, \emph{e.g.}, $e_K^p\approx 1 \cdot 10^{-2}$, the number of integration points $N$ using the optimized integration scheme is reduced substantially (in this case from 303 to 83, \emph{i.e.}, a factor of almost 4). Even substantially bigger gains are observed in three-dimensional cases \cite{divi_error-estimate-based_2020}.

From the quadrature updating patterns that emerge from the sub-cell marking strategy it is observed that, as a general trend, integration orders are increased on a per-octree-level basis. This is explained by the fact that the indicators scale with the volume of the sub-cells. Based on this observation it was anticipated that an octree-level marking strategy could be very efficient, in the sense that it would yield a similar quadrature update pattern as the sub-cell marking, but that it would need fewer iterations by virtue of marking a larger number of sub-cells per step. Fig.~\ref{fig:errorvspointsb} compares the marking strategies, conveying that the octree-level marking indeed closely follows the sub-cell marking.

Although the computational effort of the quadrature optimization algorithm is worthwhile when one wants to re-use a quadrature rule multiple times, considerable computational effort is involved. In addition, one has to set up a suitable code to determine the optimal distributions for arbitrarily cut elements. Considering this, one may not be interested in obtaining the optimized distributions of the points, but may instead want a simple rule of thumb to select the quadrature on a cut element; see, \emph{e.g.}, Refs.~\cite{abedian_finite_2013,abedian_equivalent_2019}. Using our quadrature optimization algorithm, in Ref.~\cite{divi_error-estimate-based_2020} we studied the effectivity of rules of thumb in which the order of integration is lowered with the octree depth. Although the rule-of-thumb schemes are, as expected, outperformed by the optimized schemes, they generally do provide an essential improvement in accuracy per integration point compared to equal-order integration. This observed behavior is explained by the fact that the rules of thumb qualitatively match the results of the optimization procedure.

\section{Adaptive THB-spline refinement}
\label{sec:adaptivesplines}

To leverage the flexibility of the immersed simulation paradigm with respect to refining the mesh independent of the geometry, an automated mesh adaptivity strategy is required. Various adaptivity strategies have been considered in the context of immersed methods, an overview of which is presented in, \emph{e.g.}, Ref.~\cite{divi_residual-based_2022}. These refinement strategies can be categorized as either feature-based methods (refinements are based on, \emph{e.g.}, sharp gradients in the solution or high-curvature boundary regions) or methods based on error estimates (\emph{e.g.}, residual-based or goal-oriented methods).

To develop a generic adaptive procedure for scan-based analyses, we have constructed a residual-based \emph{a posteriori} error estimator. In our isogeometric analysis approach we employ truncated hierarchical B-splines (THB-splines) \cite{giannelli_thb-splines_2012,van_brummelen_adaptive_2021} to locally refine the (volumetric) background mesh.

\subsection{Residual-based error estimation}
On account of the immersed boundary terms in the formulation (Section~\ref{sec:formulation}), it is not well-posed in the infinite dimensional setting. Upon appropriate selection of the stabilization parameters, the (mixed) Galerkin formulation of the Stokes problem is well-posed with respect to the mesh-dependent norm (see Ref.~\cite{divi_residual-based_2022} for details)
\begin{equation} \label{eq:energynormstokes}
	\ltrivert  v^h \rtrivert^2 = \ltrivert (\boldsymbol{v}^h,q^h) \rtrivert^2 =  \ltrivert \vvh \rtrivert^{2}_u + \ltrivert \qh\rtrivert^{2}_p,
\end{equation}
with 
\begin{subequations}	
	\begin{align}
		\ltrivert \vvh \rtrivert^{2}_u &:= \lVert \mu^{\frac{1}{2}} \nabla^s \boldsymbol{v}^{h} \rVert_{L^{2}(\mathcal{T})}^{2} +  \lVert \beta^{-\frac{1}{2}} h^{\frac{1}{2}}  \mu^{\frac{1}{2}} \partial_{n} \boldsymbol{v}^{h} \rVert_{L^{2}(\partial \Omega_{D})}^{2}  \nonumber \\ 
		&\phantom{:=}+ \lVert \beta^{\frac{1}{2}} h^{-\frac{1}{2}}  \mu^{\frac{1}{2}}\boldsymbol{v}^{h} \rVert_{L^{2}(\partial \Omega_{D})}^{2}+ \sum \limits_{F \in \mathcal{F}_{\rm ghost}} \lVert \gamma_g^{\frac{1}{2}} h_F^{k-\frac{1}{2}}  \mu^{\frac{1}{2}} \llbracket \partial_{n}^{k} \boldsymbol{v}^h \rrbracket \rVert_{L^{2}(F)}^{2}, \label{eq:velocitynorm} \\
		\ltrivert \qh \rtrivert^{2}_p &:=  \lVert \mu^{-\frac{1}{2}} q^{h} \rVert_{L^{2}(\mathcal{T})}^{2} + \sum \limits_{F \in \mathcal{F}_{\rm skeleton}} \lVert \gamma_{s}^{\frac{1}{2}} h_F^{k+\frac{1}{2}}  \mu^{-\frac{1}{2}} \llbracket \partial_{n}^{k} q^h \rrbracket \rVert_{L^{2}(F)}^{2}. \label{eq:pressurenorm}
	\end{align}
\end{subequations}
We refer to this mesh-dependent norm as the \emph{energy norm} and use it to construct an \emph{a posterior} error estimator for the discretization error, $u - u^h$. 

Since, in the considered immersed setting, stability can only be shown in the discrete setting, we define the solution error with respect to the solution in the order-elevated space $\widehat{u}^h \in \widehat{V}^h$. The space $\widehat{V}^h$ is defined on the same mesh and with the same regularity as the space $V^h$, but with the order of the basis elevated in such a way that $\widehat{V}^h \supset V^h$. It is then assumed that
\begin{equation}
  \ltrivert \widehat{u}^h - \uh \rtrivert \approx \ltrivert u - \uh \rtrivert.
  \label{eq:errordefinition}
\end{equation}
We note that additional stabilization terms are required to retain stability in the order-elevated space. In principle this means that the operators \eqref{eq:operators} need to be augmented, but we assume that for the solution in the order-elevated space these terms are negligible. Similar assumptions, referred to as saturation assumptions, have been considered in, \emph{e.g.}, Refs.~\cite{becker_finite_2003,juntunen_nitsches_2009,chouly_residual-based_2018}. Note that the refined space is only used to provide a proper functional setting for the error estimator and that it is not required to perform computations in this space.

To construct an estimator for the error \eqref{eq:errordefinition}, it can be bound from above by
\begin{align}
 \ltrivert \widehat{u}^h - \uh \rtrivert \lesssim \sup_{\widehat{v}^h \in \widehat{V}^h \setminus \{0\}} \frac{a^h(\widehat{u}^h-u^h,\widehat{v}^h)}{\ltrivert \widehat{v}^h \rtrivert} = \sup_{\widehat{v}^h\in \widehat{V}^h \setminus \{0\}} \frac{r^h(\widehat{v}^h)}{\ltrivert \widehat{v}^h \rtrivert},
 \label{eq:boundonenergy}
\end{align}
where the aggregate residual (\emph{i.e.}, the combined velocity-pressure residual) is defined as
\begin{equation}
\begin{aligned}
r^h(\widehat{v}^h) &:= r^h(u^h)(\widehat{v}^h) :=f^h(\widehat{v}^h) -a^h(u^h,\widehat{v}^h),
\label{eq:residualdefinition}
\end{aligned}
\end{equation}
with aggregate operators $a^h$ and $f^h$ (see Ref.~\cite{divi_residual-based_2022} for details).

We propose an error estimator pertaining to the background mesh, $\mesh$, which bounds the error in the energy norm \eqref{eq:boundonenergy} as
\begin{align}
 \mathcal{E} = \sqrt{\sum \limits_{\element \in \mesh} \eta_\element^2} \gtrsim \sup_{\widehat{v}^h\in \widehat{V}^h \setminus \{0\}} \frac{r^h(\widehat{v}^h)}{\ltrivert \widehat{v}^h \rtrivert} \gtrsim \ltrivert \widehat{u}^h - \uh \rtrivert,
 \label{eq:estimator}
\end{align}
where the element-wise error indicators, $\eta_\element$, will serve to guide an adaptive refinement procedure. 

To derive the error indicators, the residual \eqref{eq:residualdefinition} is considered with the operators defined as in \eqref{eq:operators}. Following the derivation of Ref.~\cite{divi_residual-based_2022}, the indicators are defined as
\begin{equation}
\begin{aligned}
	\eta_K^2 &=  \mu^{-1} h_K^2 \| \rintu \|_{L^2(\element \cap \domain)}^2
	+ \mu \| \rintp \|_{L^2(\element\cap \domain)}^2  \\
	 &\phantom{=} + \mu^{-1} h_K \| \rneumannu \|_{L^2({ \element \cap \neumannboundary})}^2 + 9 \mu h_K^{-1}  \| \rnitscheu \|_{L^2({\element \cap \dirichletboundary})}^2   \\
	&\phantom{=} + \sum_{\face \in \mathcal{F}_{\rm skeleton}} \mu^{-1} h_K \| \rjumpu \|_{L^2({\partial \element \cap F})}^2  +  \mu \beta^2 \meshsize_{\element}^{-1} \| \rnitscheu \|_{L^2(\element \cap \partial \domain_D)}^2 
	\Bigg. \\
	&\phantom{=} + \sum \limits_{\face \in \ghost^\element}  \mu \gamma_g^2 \meshsize_{\face}^{2\order-1}  \| \rghostu \|_{L^2(\partial \element \cap \face)}^2 \\
	&\phantom{=} + \sum \limits_{\face \in \skeleton^\element}  \mu^{-1} \gamma_s^2 \meshsize_{\face}^{2\order+1}  \| \rskeleton \|_{L^2(\partial \element \cap \face)}^2 ,
\end{aligned}
\label{eq:errorindicator}
\end{equation}
where
\begin{subequations}
\begin{align}
 \rintu &:= \bodystokes + \nabla \cdot \left( 2 \mu \nabla^s \uuh \right) - \nabla \ph,\\
 \rintp &:= \nabla \cdot \uuh,\\
 \rneumannu &:= \neumanndatastokes - \left( 2 \mu \nabla^s \uuh \right) \nn + \ph \nn ,\\
 \rnitscheu &:= \dirichletdatastokes - \uuh ,\\
 \rjumpu &:= \tfrac{1}{2}  \llbracket \left( 2 \mu \nabla^s \uuh \right) \nn \rrbracket,\\
 \rghostu &:= \tfrac{1}{2}  \llbracket \gradn^k \uuh \rrbracket,\\
 \rskeleton &:= \tfrac{1}{2} \llbracket \gradn^k \ph \rrbracket.
\end{align}%
\end{subequations}
The error indicator \eqref{eq:errorindicator} reflects that the total element error for all elements that do not intersect the boundary of the domain is composed of the interior residuals associated with the momentum balance and mass balance, and the residual terms associated with the derivative jumps on the skeleton mesh. For elements that intersect the Neumann boundary, additional error contributions are obtained from the Neumann residual and the ghost penalty residual, while additional Nitsche-related contributions appear for elements intersecting the Dirichlet boundary.

\subsection{Mesh adaptivity algorithm}
The residual-based error estimator \eqref{eq:estimator} is used in an iterative mesh refinement procedure, which is summarized in Alg.~\ref{alg:adaptivity}. The procedure takes the stabilized immersed isogeometric \textit{model} as outlined in Section~\ref{sec:formulation} as input, as well as an initial \textit{mesh} and stopping \textit{criterion}. Once the stopping \textit{criterion} is met, the algorithm returns the optimized \textit{mesh} and the corresponding \textit{solution}.

\begin{algorithm}
\caption{Function to perform residual-based error estimation and adaptivity}\label{alg:adaptivity}
    \hspace*{\algorithmicindent} \textbf{Input:} immersed isogeometric \textit{model}, initial background \textit{mesh}, stopping \textit{criterion}\\
    \hspace*{\algorithmicindent} \textbf{Output:} optimized \textit{mesh} and \textit{solution}
\begin{algorithmic}[1]
\Function{error\_estimation\_and\_adaptivity}{\textit{model}, \textit{mesh}, \textit{criterion}}
\While{\textbf{not} \textit{criterion}}\label{alg:adaptivity:whilecriterion}
\Comment{Adapt the mesh until the stopping criterion is met}
    \State \textit{solution} = \Call{solve}{\textit{model}, \textit{mesh}}\label{alg:adaptivity:solve}\Comment{Immersed IGA formulation of Section~\ref{sec:formulation}}
    \State \textit{indicators} = \Call{initialize\_indicators}{\textit{mesh}}
    \For{\textit{element}, \textit{indicator} \textbf{in} \Call{zip}{\textit{mesh}, \textit{indicators}}}
    \Comment{Iterate over the elements}\label{alg:adaptivity:elementloop}
        \State \textit{indicator} = \Call{get\_element\_indicator}{\textit{element}, \textit{model}, \textit{solution}}\label{alg:adaptivity:indicators}
    \EndFor
    \State \textit{marking} = \Call{mark\_elements}{\textit{indicators}}
    \Comment{D\"orfler marking}\label{alg:adaptivity:marking}
    \State \textit{marking} = \Call{ensure\_refinement}{\textit{marking}, \textit{mesh}}\label{alg:adaptivity:completion}
    \State \textit{mesh} = \Call{update\_mesh}{\textit{mesh}, \textit{marking}}\label{alg:adaptivity:update}
\EndWhile
\State \textbf{return} \textit{mesh}, \textit{solution}\label{alg:adaptivity:return}
\EndFunction
\end{algorithmic}
\end{algorithm}

For each step of the adaptivity procedure, for the given mesh the \textit{solution} of the Galerkin problem \eqref{eq:weak_stokes} is computed (L\ref{alg:adaptivity:solve}). For each element (L\ref{alg:adaptivity:elementloop}), the error \textit{indicator} \eqref{eq:errorindicator} is then evaluated (L\ref{alg:adaptivity:indicators}). D\"orfler marking \cite{dorfler_convergent_1996} -- targeting reduction of the estimator \eqref{eq:estimator} by a fixed fraction -- is used to select elements for refinement (L\ref{alg:adaptivity:marking}). For THB-splines, refining elements does not necessarily result in a refinement of the approximation space \cite{kuru_goal-adaptive_2014,van_brummelen_adaptive_2021}. To ensure that the approximation space is refined, an additional refinement mask  is applied to update the element \textit{marking} (L\ref{alg:adaptivity:completion}).

In our implementation the geometry approximation is not altered during mesh refinement. A consequence of this implementation choice is that an element can only be refined up to the octree depth. Elements requiring refinement beyond this depth are discarded from the \textit{marking} list, and the adaptive refinement procedure is stopped if there are no more elements that can be refined. We refer the reader to Ref.~\cite{divi_residual-based_2022} for details.

\subsection{Mesh adaptivity results}
\begin{figure}
	\centering
	\begin{subfigure}[b]{0.5\textwidth}
		\centering
		\includegraphics[width=0.8\textwidth]{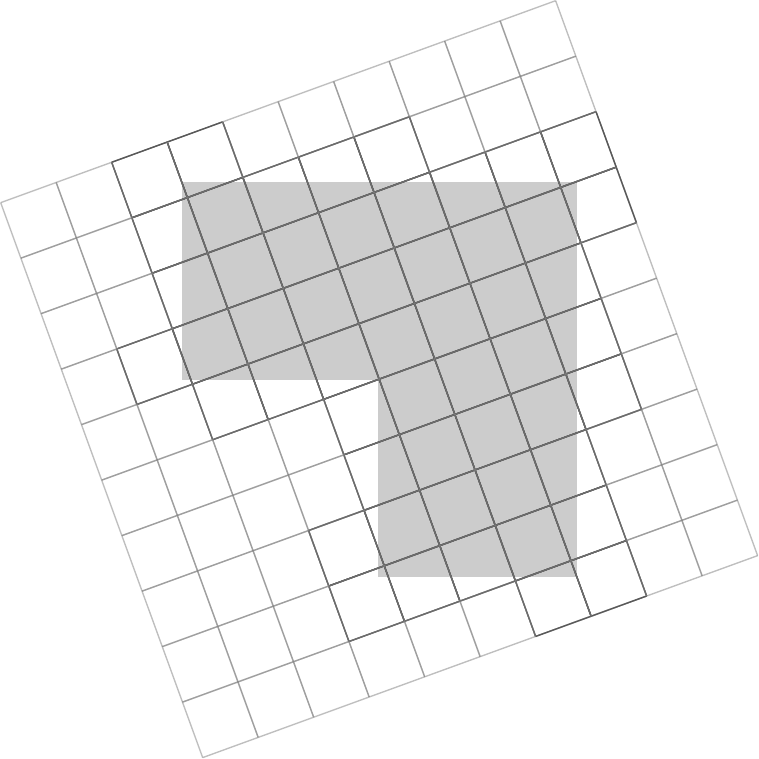}
		\caption{Initial mesh}\label{fig:Lshape_stokes_adap_mesha}
	\end{subfigure}\hfill%
	\begin{subfigure}[b]{0.5\textwidth}
		\centering
		\includegraphics[width=0.8\textwidth]{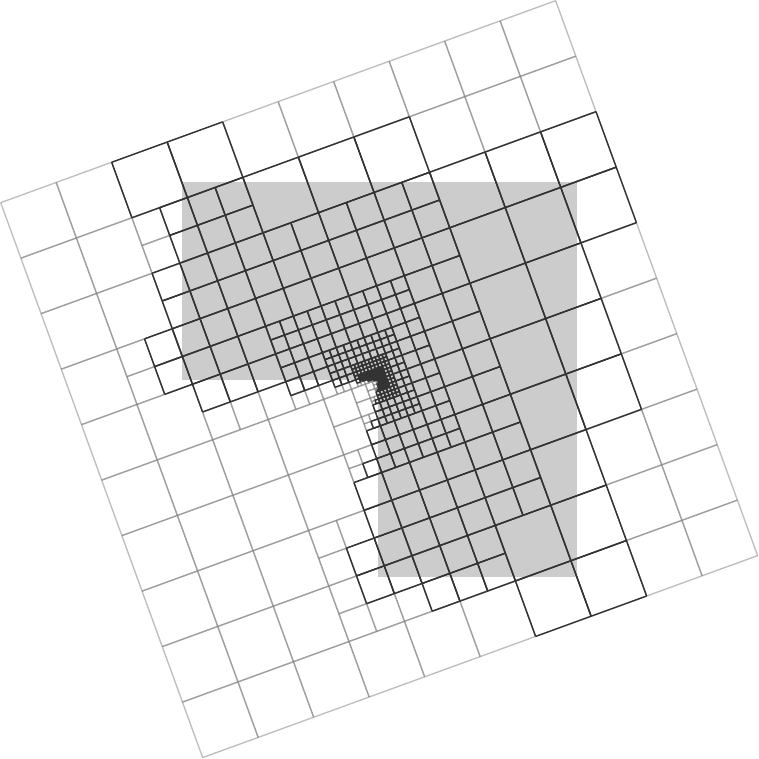}
		\caption{Step $5$}\label{fig:Lshape_stokes_adap_meshb}
	\end{subfigure}
	\caption{Evolution of the mesh using the adaptive refinement procedure for the Stokes problem on a re-entrant corner domain using $\order = 2$.}
	\label{fig:Lshape_stokes_adap_mesh}
\end{figure}

Before considering the application of the developed residual-based error estimation and adaptivity procedure in the context of scan-based analysis in Section~\ref{sec:results}, we here first reproduce a benchmark case from Ref.~\cite{divi_residual-based_2022}. We consider the Stokes problem \eqref{eq:stokesequations} on a re-entrant corner domain (Fig.~\ref{fig:Lshape_stokes_adap_mesha}) with mixed Dirichlet and Neumann boundaries. The method of manufactured solutions is considered with the (weakly singular) exact solution taken from Ref.~\cite{verfurth_review_1996}. We refer to Ref.~\cite{divi_residual-based_2022} and references therein for a full specification of the benchmark.

\begin{figure}
	\centering
	\begin{subfigure}[b]{0.5\textwidth}
		\centering
		\includegraphics[width=\textwidth]{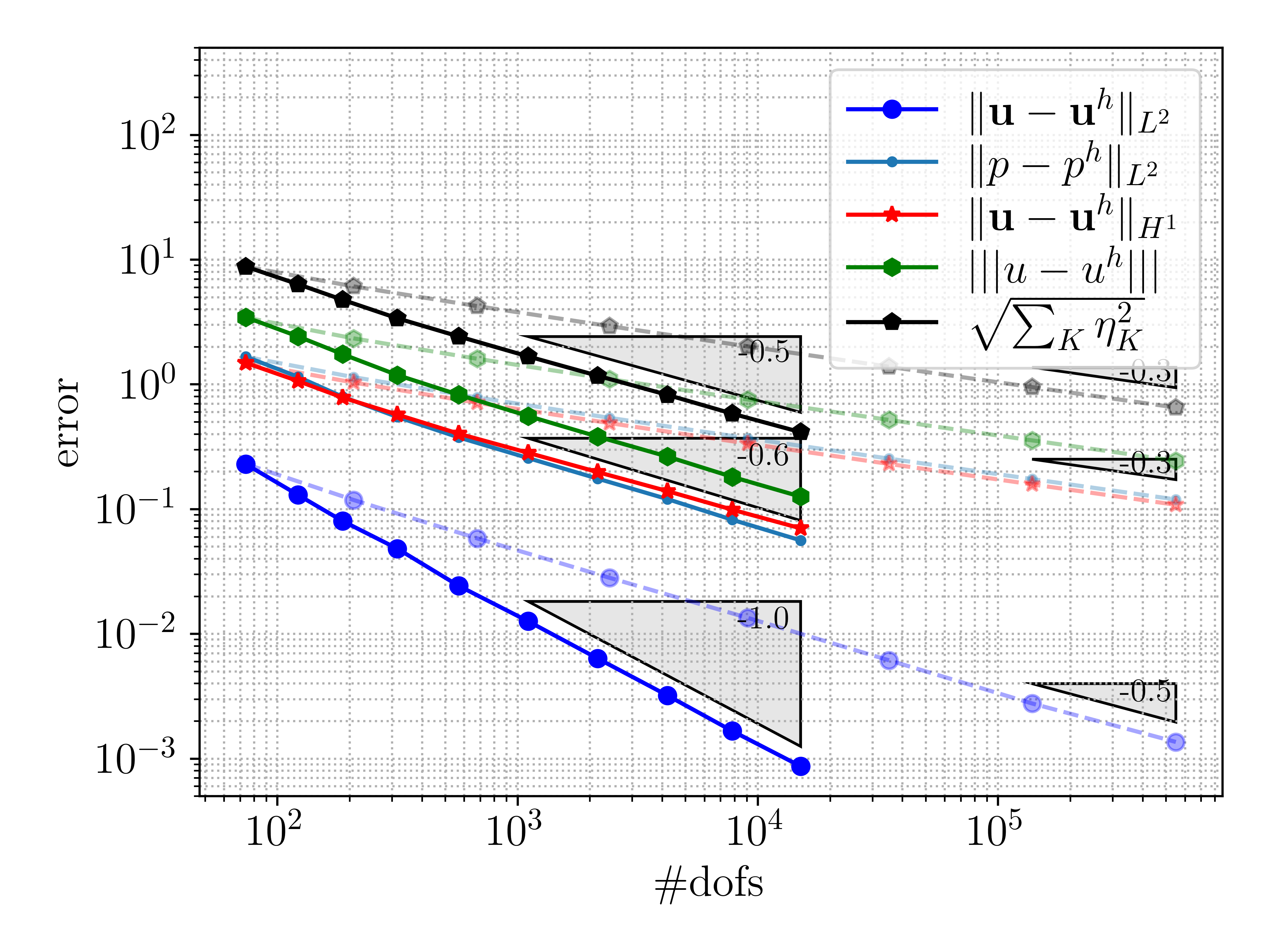}
		\caption{$\order = 1$}
		\label{fig:Lshape_stokes_p1}
	\end{subfigure}\hfill
	\begin{subfigure}[b]{0.5\textwidth}
		\centering
		\includegraphics[width=\textwidth]{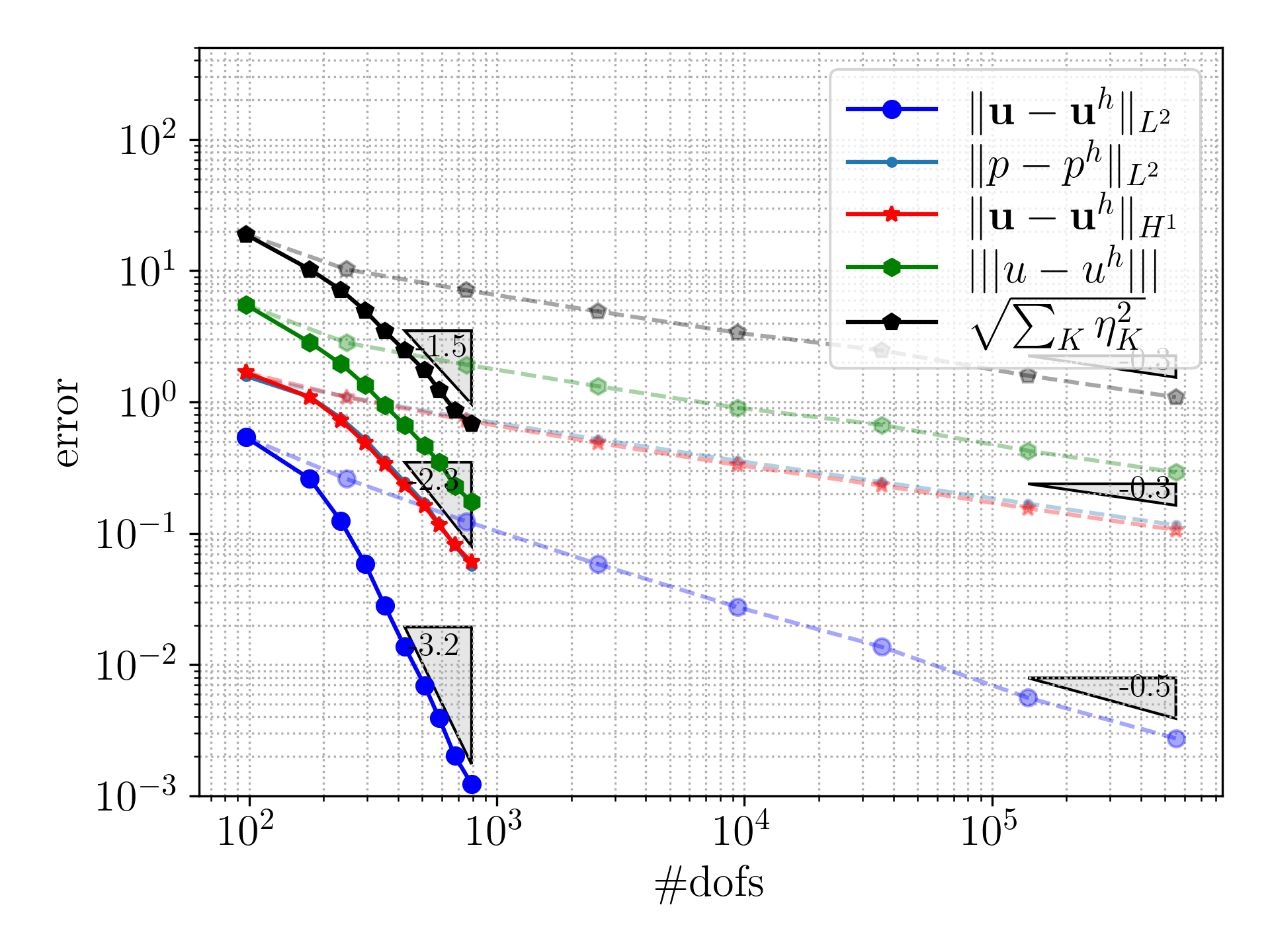}
		\caption{$\order = 2$}
		\label{fig:Lshape_stokes_p2}
	\end{subfigure}
		\caption{Error convergence results for the Stokes problem on a re-entrant corner domain under residual-based adaptive refinement (solid) and uniform refinement (dashed) for linear ($\order =1$) and quadratic ($\order =2$) basis functions.}
	\label{fig:Lshape_stokes_conv}
\end{figure}

Fig.~\ref{fig:Lshape_stokes_conv} displays the error convergence results obtained using uniform and adaptive refinements, for both linear and quadratic THB-splines. The convergence rate when uniform refinements are considered is suboptimal, limited by the weak singularity at the re-entrant corner. Using adaptive mesh refinement results in a recovery of the optimal rates in the case of linear basis functions, with even higher rates observed for the quadratic splines on account of the highly-focused refinements resulting from the residual-based error estimator as observed in Fig.~\ref{fig:Lshape_stokes_adap_meshb}.

\section{Scan-based flow simulations}
\label{sec:results}

To demonstrate the scan-based analysis workflow reviewed in this work, we consider the blood flow (viscosity $\mu=4\,{\rm mPa \cdot s}$) through the patient-specific ${\rm \mu CT}$-based carotid artery introduced in Section~\ref{sec:topology} (Fig.~\ref{fig:topology}). Neumann conditions are imposed on the inflow (bottom) and outflow (top) boundaries, with the traction on the inflow boundary corresponding to a pressure of $17.3$\,kPa ($130$\,mmHg) and a zero-traction condition on the outflow boundary. Homogeneous Dirichlet conditions are imposed along the immersed boundaries to impose a no slip condition. The presented results are based on second-order ($\order = 2$) THB-splines. For details regarding the simulation setup we refer the reader to Ref.~\cite{divi_residual-based_2022}, from which the results presented here are reproduced.

\begin{figure}
	\begin{subfigure}[b]{\textwidth}
		\centering
		\includegraphics[height=0.49\textwidth]{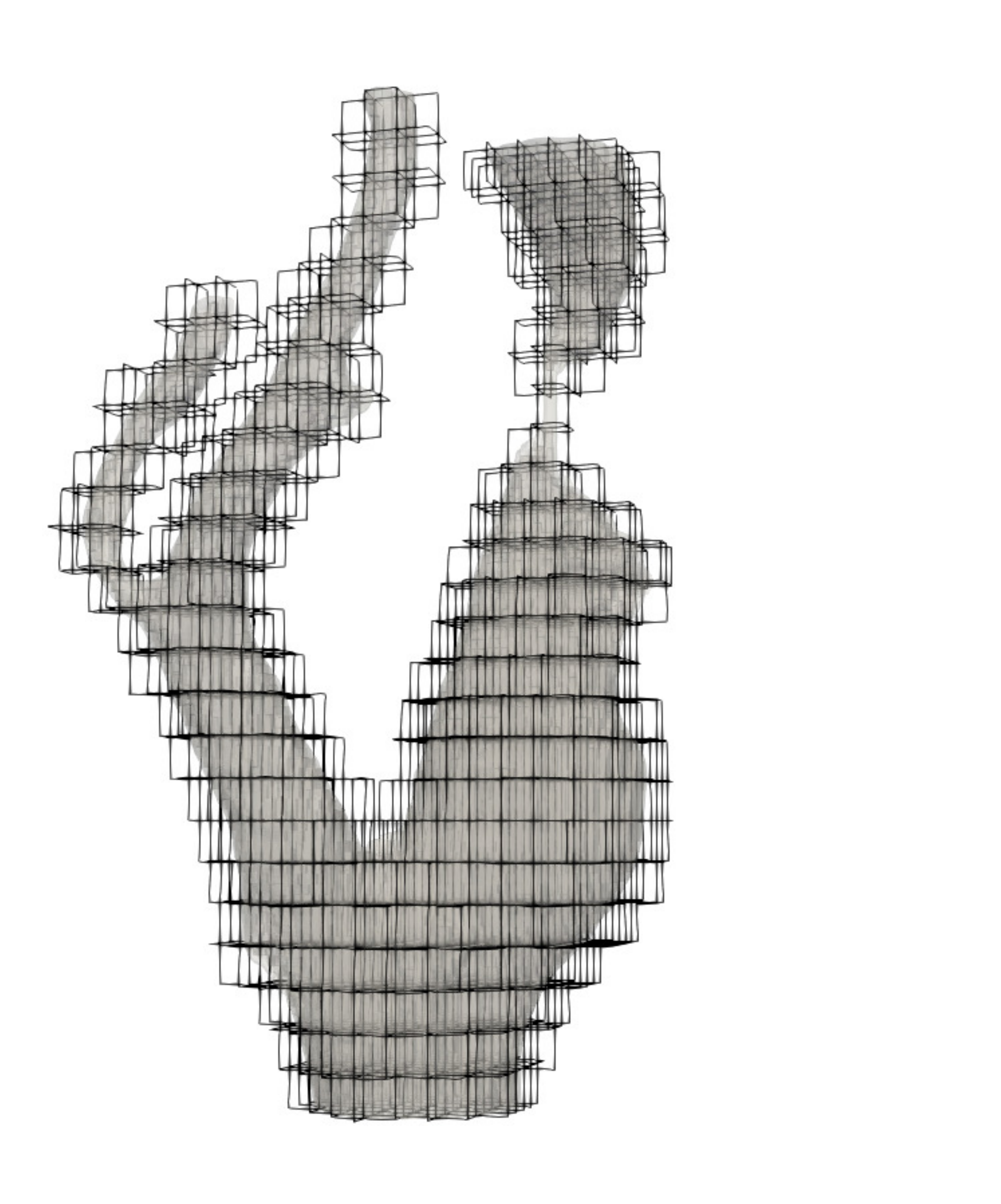}
		\includegraphics[height=0.49\textwidth]{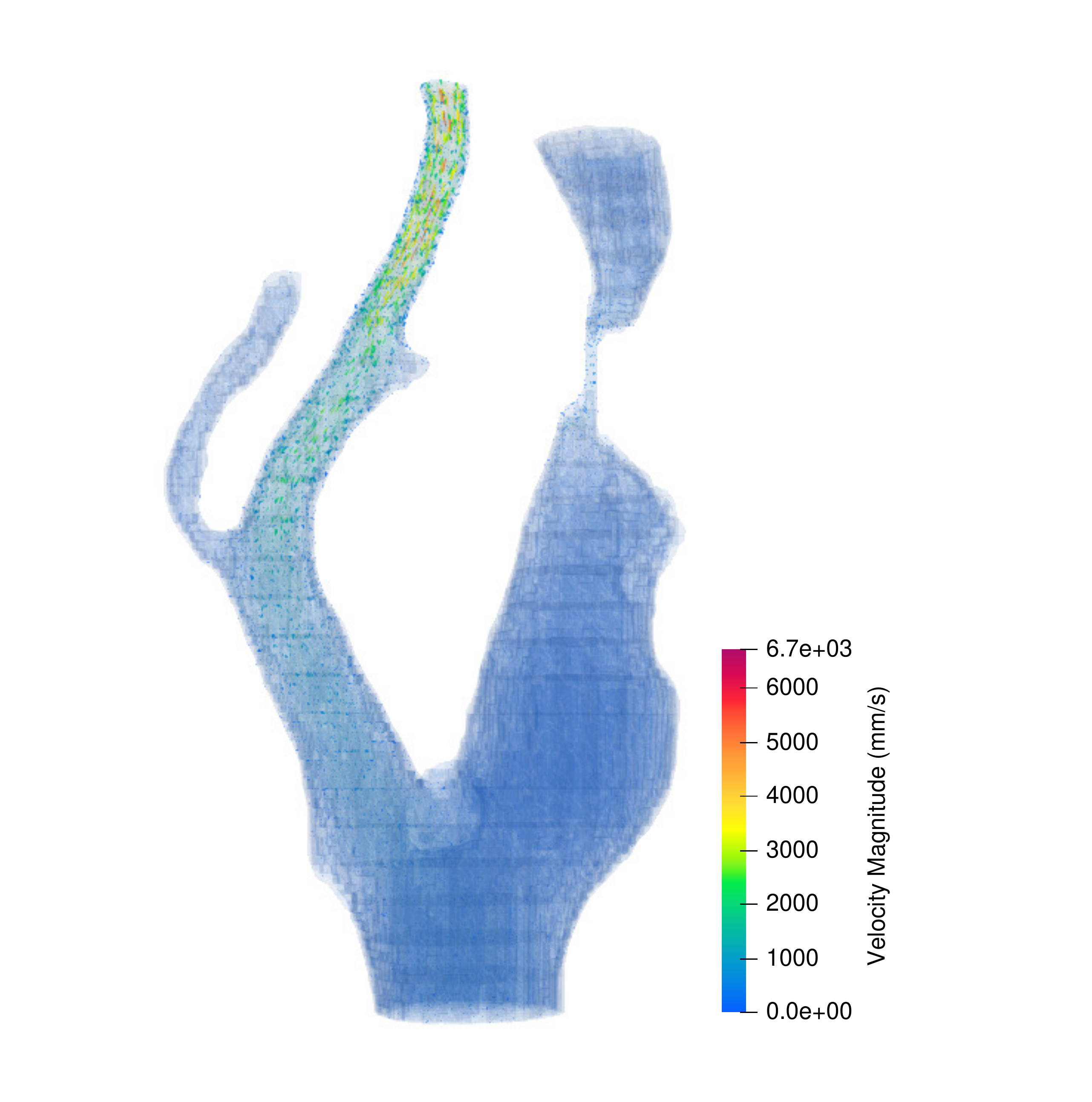}
		\caption{Initial mesh with $3158$ \#DOFs}
		\label{fig:3dstokes_velocity_step0}
	\end{subfigure}
	\begin{subfigure}[b]{\textwidth}
		\centering
		\includegraphics[height=0.49\textwidth]{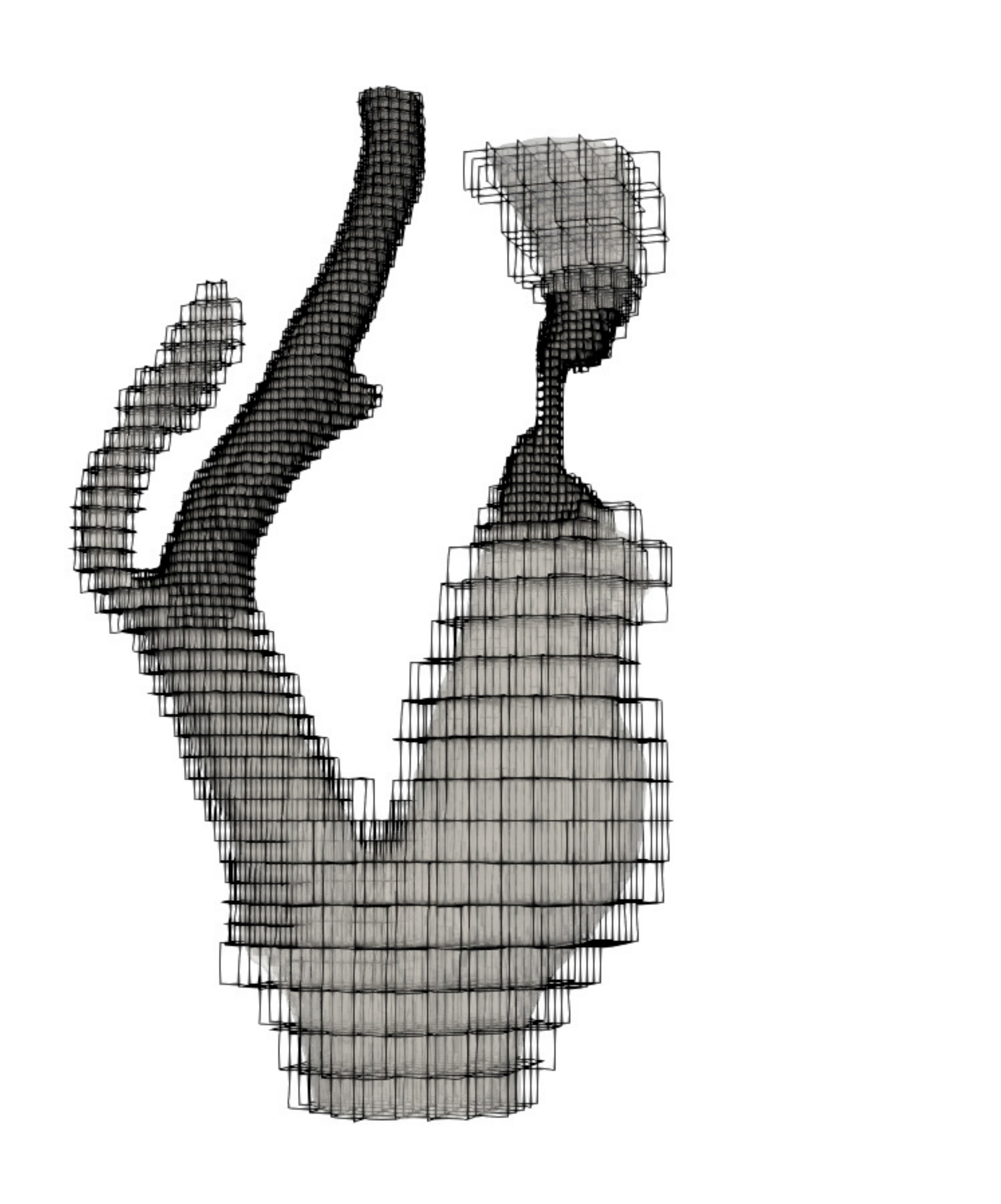}
		\includegraphics[height=0.49\textwidth]{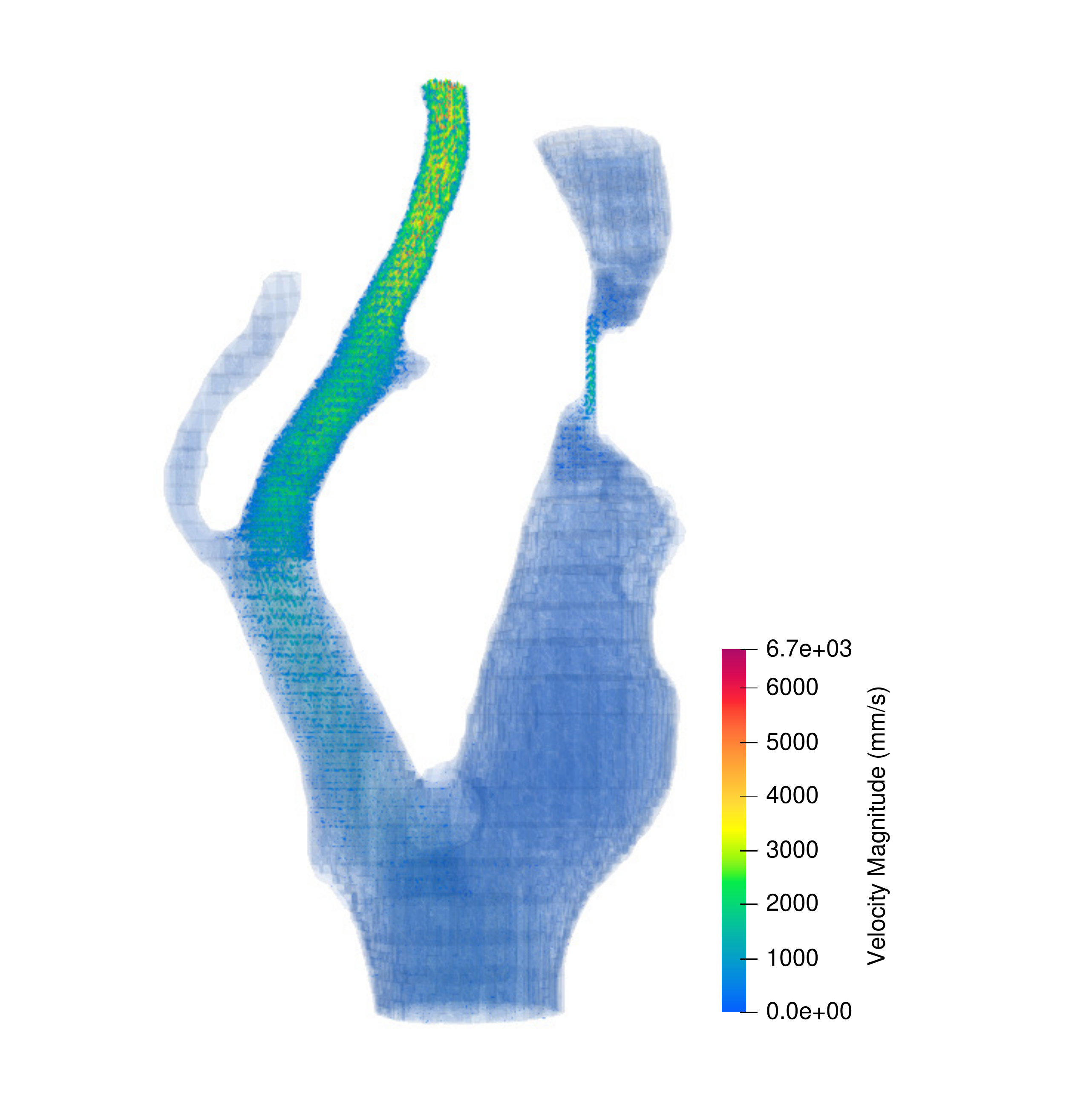}
		\caption{Step $3$ with  $12467$ \#DOFs}
		\label{fig:3dstokes_velocity_step3}
	\end{subfigure}
	\caption{Initial (\subref{fig:3dstokes_velocity_step0}) and final (\subref{fig:3dstokes_velocity_step3}) mesh (left) and velocity magnitude (right) for the patient-specific flow problem.}
	\label{fig:3dstokes_velocity}
\end{figure}

We consider an initial scan-domain mesh consisting of $24 \times 24 \times 24$ elements, with a scan size of $25.6 \times 21.1 \times 32.0\, {\rm mm}^3$. The octree depth is set to three. In this setting, after two refinements, an element is of a similar size as the voxels. The need to substantially refine beyond the voxel size is, from a practical perspective, questionable, as the dominant error in the analysis will then be related to the scan resolution and the segmentation procedure. In this sense, the constraint of not being able to refine beyond the octree depth is not a crucial problem in the considered simulations.

The initial mesh and final refinement result are shown in Fig.~\ref{fig:3dstokes_velocity}. The adaptive refinement procedure focuses on the regions where the errors are largest, \emph{i.e.}, near the stenosed section (\emph{i.e.}, the narrow region at the right artery) and at the outflow section of the left artery, such that important details of the solution are resolved. After the final refinement step, the adaptive simulation uses $12,816$ DOFs, which is substantially lower than the approximately $100,000$ DOFs that would have resulted from uniform mesh refinements up to the same level \cite{divi_topology-preserving_2022}.

\section{Concluding remarks}
\label{sec:conclusions}

In this contribution, we have reviewed the four key research contributions of our team with respect to scan-based immersed isogeometric flow analysis, \emph{viz.}: \emph{(i)} A spline-based image segmentation procedure, encompassing a voxel-data smoothing procedure, an octree-based procedure to obtain an explicit parametrization of the computational domain and its (immersed) boundary, and a topology-preservation strategy to restore smoothing-induced anomalies; \emph{(ii)} A stabilized immersed formulation for (\emph{a.o.}) Stokes flow, which ensures robustness with respect to unfavorably cut elements and enables the consideration of equal-order velocity-pressure discretizations without the loss of inf-sup stability; \emph{(iii)} An adaptive procedure to optimize the distribution of integration points over cut elements, based on Strang's first lemma; \emph{(iv)} A mesh refinement procedure based on rigorous residual-based error estimates to refine the computational mesh in places where this results in significant accuracy improvements.

An important aspect of immersed (finite element) methods is the ill-conditioning associated with small (\emph{i.e.}, with a small volume fraction) or unfavorably cut (\emph{e.g.}, sliver-like) elements. Although not reviewed in this work, over the past decade our team has contributed to solving the challenges associated with ill-conditioning. The origin of the ill-conditioning problem was studied in detail by De Prenter \emph{et al.} \cite{de_prenter_condition_2017}, which led to a scaling relation for the condition number with the smallest cut-element volume fraction. Dedicated preconditioning techniques, to be used in conjunction with iterative solvers, were developed based on the insights from this work, \emph{e.g.}, Refs.~\cite{de_prenter_condition_2017,de_prenter_preconditioning_2019,de_prenter_multigrid_2020}. We consider these (preconditioned) solver developments an important step in unlocking the potential of high-performance computing for immersed finite element methods \cite{jomo_robust_2019}. Note that the ghost- and skeleton-stabilization terms employed in the formulation in this chapter, which are primarily added to ensure well-posedness of the weak form, also resolve the conditioning problems, such that preconditioning techniques are not essential in this work. 

The innovations in computational procedures and problem formulations yield a highly robust immersed isogeometric analysis workflow for scan-based analyses. Error-controlled simulations can be performed directly based on scan data, without the need for extensive user interactions. The effectivity of the framework is not fundamentally affected by the geometric and topological complexity of the scan data, on account of the decoupling of the geometry and computational mesh in immersed methods. The robustness of the framework derives from the rigorous mathematical underpinning of the considered methods.

Further developments to the scan-based workflow are required to enable the consideration of more advanced problems/formulations, such as higher Reynolds number flows (requiring additional stabilization), fluid-structure interactions and complex fluid models. Further improvements are also required to enhance the computational performance of the developed workflow. This mainly pertains to algorithm and code optimization, which is required to apply the developed workflow to, \emph{e.g.}, larger scans, time-dependent problems and non-linear problems. Detailed recommendations for specific further developments can be found in our referenced work; see Ref.~\cite{divi_scan-based_2022} for a summary of these.

\begin{acknowledgement}
Our implementation is based on the open source finite element library Nutils \cite{zwieten_nutils_2022}. CVV and SCD acknowledge the partial support of the European Union's Horizon 2020 research and innovation programme under Grant Agreement No 101017578 (SIMCor).
\end{acknowledgement}

\printbibliography

\end{document}